\documentclass[11pt]{amsart}
\usepackage{amsxtra}
\usepackage{amssymb}
\theoremstyle{plain}

\newcommand{\cleqn}{\setcounter{equation}{0}}
\newcommand{\clth}{\setcounter{theorem}{0}}
\newcommand {\sectionnew}[1]{\section{#1}\cleqn\clth}

\newtheorem{theorem}{Theorem}[section]
\newtheorem{lemma}[theorem]{Lemma}
\newtheorem{definition-theorem}[theorem]{Definition-Theorem}
\newtheorem{proposition}[theorem]{Proposition}
\newtheorem{corollary}[theorem]{Corollary}
\newtheorem{definition}[theorem]{Definition}
\newtheorem{example}[theorem]{Example}
\newtheorem{remark}[theorem]{Remark}
\newtheorem{conjecture}[theorem]{Conjecture}
\newtheorem{notation}[theorem]{Notation}

\newcommand \bth[1] { \begin{theorem}\label{t#1} }
\newcommand \ble[1] { \begin{lemma}\label{l#1} }

\newcommand \bpr[1] { \begin{proposition}\label{p#1} }
\newcommand \bco[1] { \begin{corollary}\label{c#1} }
\newcommand \bde[1] { \begin{definition}\label{d#1}\rm }
\newcommand \bex[1] { \begin{example}\label{e#1}\rm }
\newcommand \bre[1] { \begin{remark}\label{r#1}\rm }
\newcommand \bcj[1] { \begin{conjecture}\label{j#1}\rm }

\newcommand \bnota[1] { \begin{notation}\label{n#1}\rm }
\renewcommand {\eth} { \end{theorem} }
\newcommand {\ele} { \end{lemma} }

\newcommand {\epr} { \end{proposition} }
\newcommand {\eco} { \end{corollary} }
\newcommand {\ede} { \end{definition} }
\newcommand {\eex} { \end{example} }
\newcommand {\ere} { \end{remark} }
\newcommand {\ecj} { \end{conjecture} }

\newcommand {\enota} { \end{notation} }
\newcommand \thref[1]{Theorem \ref{t#1}}
\newcommand \leref[1]{Lemma \ref{l#1}}
\newcommand \prref[1]{Proposition \ref{p#1}}
\newcommand \coref[1]{Corollary \ref{c#1}}

\newcommand \deref[1]{Definition \ref{d#1}}
\newcommand \exref[1]{Example \ref{e#1}}

\def \KK {{\mathbb K}}
\def \Zset {{\mathbb Z}}
\def \Qset {{\mathbb Q}}
\def \Nset {{\mathbb N}}

\def \AA {{\mathcal{A}}}
\def \NN {{\mathcal{N}}}

\def \TT {{\mathcal{T}}}

\def \CC {{\mathcal{C}}}

\def \OO {{\mathcal{O}}}

\def \UU {{\mathcal{U}}}

\def \SS {{\mathcal{S}}}

\def \TT {{\mathcal{T}}} 

\def \pb {{\bf{p}}}
\def \qb {{\bf{q}}}
\def \de {\delta}
\def \al {\alpha}
\def \be {\beta}
\def \vpi {\varpi}
\def \la {\lambda}

\def \Om {\Omega}
\def \ga {\gamma}
\def \de {\delta}

\def \sig {\sigma}

\def \vph {\varphi}

\def \sig{\sigma}

\def \lab {\boldsymbol{\lambda}}

\def \mt  {\mapsto}

\def \rcor {\rangle}
\def \lcor {\langle}

\def \ol {\overline}

\def \id { {\mathrm{id}} }

\DeclareMathOperator \rank { \operatorname{rank} }
\DeclareMathOperator \UAut { \operatorname{UAut} }
\DeclareMathOperator \DAut { \operatorname{DAut} }

\def \g  {\mathfrak{g}}   

\def \sl {\mathfrak{sl}}

\def \n  {\mathfrak{n}}

\def \sl {\mathfrak{sl}}

\DeclareMathOperator \Aut { {\mathrm{Aut}} }

\DeclareMathOperator \Ext { {\mathrm{Ext}} }

\DeclareMathOperator \Fract { {\mathrm{Fract}} }

\DeclareMathOperator \rk {\operatorname{rank}}
\DeclareMathOperator \rad {\operatorname{rad}}
\newcommand\kx{\KK^*}

\newcommand\HH{{\mathcal{H}}}
\newcommand\xh{X(\HH)}

\newcommand \Znn {\Zset_{\ge 0}}

\newcommand \Oq {\OO_q}
\begin{document}
\title[Unipotent and Nakayama automorphisms of quantum nilpotent algebras]
{Unipotent and Nakayama automorphisms of \\
quantum nilpotent algebras}
\author[K. R. Goodearl]{K. R. Goodearl}
\address{
Department of Mathematics \\
University of California\\
Santa Barbara, CA 93106 \\
U.S.A.
}
\email{goodearl@math.ucsb.edu}
\author[M. T. Yakimov]{M. T. Yakimov}
\thanks{The research of K.R.G. was partially supported by NSF grant DMS-0800948, 
and that of M.T.Y. by NSF grant DMS-1001632 and and DMS-1303036.}
\address{
Department of Mathematics \\
Louisiana State University \\
Baton Rouge, LA 70803 \\
U.S.A.
}
\email{yakimov@math.lsu.edu}
\date{}
\keywords{Iterated skew polynomial algebras, iterated Ore extensions, CGL extensions, noncommutative UFDs, unipotent automorphisms, twisted Calabi-Yau algebras, Nakayama automorphisms}
\subjclass[2010]{Primary 16T20; Secondary 16S36, 16W20, 17B37, 20G42}

\begin{abstract} 
Automorphisms of algebras $R$ from a very large axiomatic class of quantum nilpotent algebras are studied using techniques 
from noncommutative unique factorization domains and quantum cluster algebras. 
First, the Nakayama automorphism of $R$ (associated to its structure
as a twisted Calabi-Yau algebra) is determined and shown to be given by
conjugation by a normal element, namely, the product of the
homogeneous prime elements of $R$ (there are finitely many up to
associates).
Second, in the case when $R$ is connected graded, the unipotent automorphisms of $R$ are classified up to minor exceptions. This theorem is a far reaching extension of the classification results \cite{Y-LL,Y-AD}
previously used to settle the Andruskiewitsch--Dumas and Launois--Lenagan conjectures. 
The result on unipotent automorphisms has a wide range of applications to the determination 
of the full automorphisms groups of the connected graded algebras in the family. This is 
illustrated by a uniform treatment of the automorphism groups of the generic algebras of quantum matrices
of both rectangular and square shape \cite{LL1,Y-LL}.
\end{abstract}

\maketitle

\sectionnew{Introduction}
\label{intro}

This paper is devoted to a study of automorphisms of \emph{quantum nilpotent algebras}, a large, axiomatically defined class of algebras. The algebras in this class are known under the name Cauchon--Goodearl--Letzter extensions and consist of iterated 
skew polynomial rings satisfying certain common properties for algebras appearing in the area of quantum groups.
The class contains the quantized coordinate rings of the Schubert cells for all simple algebraic groups, multiparameter quantized coordinate rings of many algebraic varieties,
quantized Weyl algebras, and related algebras. The quantized coordinate rings of all double Bruhat cells are localizations of 
special algebras in the class.  

Extending the results of \cite{AC,LL1,Y-LL,Y-AD}, we prove that all of these algebras are relatively rigid in terms of symmetry, in the sense that they have far fewer automorphisms than their classical counterparts. This allows strong control, even exact descriptions in many cases, of their automorphism groups. We pursue this theme in two directions. First, results of Liu, Wang, and Wu \cite{LWW} imply that any quantum nilpotent algebra $R$ is a twisted Calabi-Yau algebra. In particular, $R$ thus has a special associated automorphism, its \emph{Nakayama automorphism}, which controls twists appearing in the cohomology of $R$. At the same time, 
all algebras $R$ in the class that we consider are equivariant noncommutative unique factorization domains \cite{LLR} in the sense of Chatters \cite{Cha}.
We develop a formula for the Nakayama automorphism $\nu$ of $R$, and show that $\nu$ is given by commutation with a special normal element. Specifically, if $u_1,\dots,u_n$ is a complete list of the homogeneous prime elements of $R$ up to scalar multiples, then $a (u_1\cdots u_n) = (u_1\cdots u_n) \nu(a)$ for all $a \in R$. (Here homogeneity is with respect to the grading of $R$ arising from an associated torus action.) It was an open problem to understand what is the role of the special element of 
the equivariant UFD $R$ that equals the product of all (finitely many up to associates) homogeneous prime elements of $R$. The first main result in the paper answers this: conjugation by this special element is the Nakayama automorphism of $R$.  

In a second direction, we obtain very general rigidity results for the connected graded algebras $R$ in the abovementioned axiomatic class. This is done by combining the quantum cluster algebra structures 
that we constructed in \cite{GY,GY2} with the rigidity of quantum tori theorem of \cite{Y-AD}. 
The quantum clusters of $R$ constructed in \cite{GY,GY2} provide a huge supply of embeddings 
$\AA \subseteq R \subset \TT$ where $\AA$ is a quantum affine space algebra and $\TT$ 
is the corresponding quantum torus. 
This allows for strong control of the \emph{unipotent} automorphisms of $R$ relative to a nonnegative grading on $R$, those being automorphisms $\psi$ such that for any homogeneous element $x \in R$ of degree $d$, the difference $\psi(x) - x$ is supported in degrees greater than $d$. Such a $\psi$ induces \cite{Y-LL} a continuous automorphism of the completion of any quantum torus $\TT$ as above, to which a general rigidity theorem proved in \cite{Y-AD} applies. We combine this rigidity with the large supply of quantum clusters 
in \cite{GY,GY2} and a general theorem for separation of variables from \cite{GY}. With this combination of methods and
the noncommutative UFD property of $R$, we show here that the unipotent automorphisms of a quantum nilpotent algebra $R$ 
have a very restricted form, which is a very general improvement of the earlier results in that direction \cite{Y-LL, Y-AD}
that were used in proving the Andruskiewitsch--Dumas and Launois--Lenagan conjectures.
Our theorem essentially classifies the unipotent automorphisms of all connected graded algebras 
in the class, up to the presence of certain types of torsion in the scalars involved in the algebras. In a variety of cases the full automorphism group $\Aut(R)$ can be completely  determined as an application of this result. We illustrate this by presenting, among other examples, a new route to the determination of the automorphism groups of generic quantum matrix algebras \cite{LL1, Y-LL} of both rectangular and square shape, 
in particular giving a second proof of the conjecture in \cite{LL1}.

In a recent paper \cite{CPWZ}, Ceken, Palmieri, Wang, and Zhang classified 
the automorphism groups of certain PI algebras using discriminants. 
Their methods apply to quantum affine spaces at roots of unity but not to general quantum matrix algebras 
at roots of unity. It is an interesting problem whether the methods of quantum cluster algebras and 
rigidity of quantum tori can be applied in conjunction with the methods of \cite{CPWZ} to treat the automorphism 
groups of the specializations of all algebras in this paper to roots of unity.

We finish the introduction by describing the class of quantum nilpotent algebras 
that we address. These algebras are iterated skew polynomial extensions
\begin{equation} 
\label{itOre.intro}
R := \KK[x_1][x_2; \sig_2, \delta_2] \cdots [x_N; \sig_N, \delta_N] 
\end{equation}
over a base field $\KK$, equipped with rational actions by tori $\HH$ of automorphisms which cover the $\sig_k$ in a suitably generic fashion, and such that the skew derivations $\de_k$ are locally nilpotent. They have been baptized \emph{CGL extensions} in \cite{LLR}; see \deref{CGL} for the precise details. We consider the class of CGL extensions to be the best current definition of quantum 
nilpotent algebras from a ring theoretic perspective. All important CGL extensions that we are aware of are \emph{symmetric} in the sense that they possess CGL extension presentations with the generators in both forward and reverse orders, that is, both \eqref{itOre.intro} and 
$$
R = \KK[x_N] [x_{N-1}; \sig^*_{N-1}, \de^*_{N-1}] \cdots [x_1; \sig^*_1, \de^*_1].
$$
The results outlined above apply to the class of symmetric CGL extensions satisfying a mild additional assumption on the scalars that appear.
\medskip

Throughout, fix a base field $\KK$. All automorphisms are assumed to be $\KK$-algebra automorphisms, and all skew derivations are assumed to be $\KK$-linear. We also assume that in all Ore extensions (skew polynomial rings) $B[x;\sig,\de]$, the map $\sig$ is an automorphism. Recall that $B[x;\sig,\de]$ denotes a ring generated by a unital subring $B$ and an element $x$ satisfying $xs = \sig(s) x + \de(s)$ for all $s \in S$, where $\sig$ is an automorphism of $B$ and $\de$ is a {\rm(}left\/{\rm)} $\sig$-derivation of $B$.

We will denote $[j,k] := \{ n \in \Zset \mid j \le n \le k \}$ for $j,k \in \Zset$. In particular, $[j,k] = \varnothing$ if $j \nleq k$.
\medskip
\\
\noindent
{\bf Acknowledgements.} 
We would like to thank MSRI for its hospitality during the programs 
in ``Cluster Algebras'' (2012) and ``Noncommutative Algebraic Geometry and Representation 
Theory'' (2013) when parts of this project were carried out.
\sectionnew{Symmetric CGL extensions}
\label{sectCGL}

In this section, we give some background on $\HH$-UFDs and CGL extensions, including some known results, and then establish a few additional results that will be needed in later sections.

\subsection{$\HH$-UFDs}  \label{HUFD}

Recall that a {\em{prime element}} of a domain $R$ is any nonzero normal element $p\in R$ (\emph{normality} meaning that $Rp = pR$)
such that $Rp$ is a completely prime ideal, i.e., $R/Rp$ is a domain.
Assume that in addition $R$ is a $\KK$-algebra and $\HH$ a group acting on $R$ 
by $\KK$-algebra automorphisms. An {\em{$\HH$-prime ideal}} of $R$ is any proper 
$\HH$-stable ideal $P$ of $R$ such that $(IJ\subseteq P \implies I\subseteq P$ 
or $J\subseteq P)$ for all $\HH$-stable ideals $I$ and $J$ of $R$. In general, $\HH$-prime ideals need not be prime, but they are prime in the case of CGL extensions \cite[II.2.9]{BG}.

One says that $R$ is an {\em{$\HH$-UFD}} if each nonzero $\HH$-prime ideal 
of $R$ contains a prime $\HH$-eigenvector. This is an equivariant 
version of Chatters' notion \cite{Cha} of noncommutative unique 
factorization domain given in \cite[Definition 2.7]{LLR}.

The following fact is an equivariant version of results of 
Chatters and Jordan \cite[Proposition 2.1]{Cha},
\cite[p. 24]{ChJo}, see \cite[Proposition 2.2]{GY} and
\cite[Proposition 6.18 (ii)]{Y-sqg}. 

\bpr{factorHUFD}
Let $R$ be a noetherian $\HH$-UFD. Every normal $\HH$-eigenvector in $R$ 
is either a unit or a product of prime $\HH$-eigenvectors. The factors are unique 
up to reordering and taking associates.
\epr

We shall also need the following equivariant version of \cite[Lemma 2.1]{ChJo}.
A nonzero ring $R$ equipped with an action of a group $\HH$ is said to be \emph{$\HH$-simple} provided the only $\HH$-stable ideals of $R$ are $0$ and $R$.

\ble{IcapER}
Let $R$ be a noetherian $\HH$-UFD and $E(R)$ the multiplicative subset of $R$ generated by the prime $\HH$-eigenvectors of $R$. All nonzero $\HH$-stable ideals of $R$ meet $E(R)$, and so the localization $R[E(R)^{-1}]$ is $\HH$-simple.
\ele

\begin{proof} The second conclusion is immediate from the first. To see the first, let $I$ be a nonzero $\HH$-stable ideal of $R$. Since $R$ is noetherian, $P_1 P_2 \cdots P_m \subset I$ for some prime ideals $P_j$ minimal over $I$. For each $j$, the intersection of the $\HH$-orbit of $P_j$ is an $\HH$-prime ideal $Q_j$ of $R$ such that $I \subseteq Q_j \subseteq P_j$. Each $Q_j$ contains a prime $\HH$-eigenvector $q_j$, and the product $q_1 q_2 \cdots q_m$ lies in $I$. Thus, $I \cap E(R) \ne \varnothing$, as desired. (Alternatively, suppose that $I \cap E(R) = \varnothing$, enlarge $I$ to an $\HH$-stable ideal $P$ maximal with respect to being disjoint from $E(R)$, check that $P$ is $\HH$-prime, and obtain a prime $\HH$-eigenvector in $P$, yielding a contradiction.)
\end{proof}

\subsection{CGL extensions} \label{genCGL}
Throughout the paper, we focus on iterated Ore extensions of the form
\begin{equation} 
\label{itOre}
R := \KK[x_1][x_2; \sig_2, \delta_2] \cdots [x_N; \sig_N, \delta_N] .
\end{equation}
We refer to such an algebra as an \emph{iterated Ore extension over $\KK$}, to emphasize that the initial step equals the base field $\KK$, and we call the integer $N$ the \emph{length} of the extension. For $k \in [0,N]$, we let $R_k$ denote the subalgebra of $R$ generated by $x_1, \dots, x_k$. In particular, $R_0 = \KK$ and 
$R_N = R$. Each $R_k$ is an iterated Ore extension over $\KK$, of length $k$.

\bde{CGL} An iterated Ore extension \eqref{itOre}
is called a \emph{CGL extension} 
\cite[Definition 3.1]{LLR} if it is equipped with a rational action of a $\KK$-torus $\HH$ 
by $\KK$-algebra automorphisms satisfying the following conditions:
\begin{enumerate}
\item[(i)] The elements $x_1, \ldots, x_N$ are $\HH$-eigenvectors.
\item[(ii)] For every $k \in [2,N]$, $\de_k$ is a locally nilpotent 
$\sig_k$-derivation of the algebra $R_{k-1}$. 
\item[(iii)] For every $k \in [1,N]$, there exists $h_k \in \HH$ such that 
$\sig_k = (h_k \cdot)|_{R_{k-1}}$ and the $h_k$-eigenvalue of $x_k$, to be denoted by $\la_k$, is not a root of unity.
\end{enumerate}

Conditions (i) and (iii) imply that 
$$
\sig_k(x_j) = \la_{kj} x_j \; \; \mbox{for some} \; \la_{kj} \in \kx, \; \; \forall 1 \le j < k \le N.
$$
We then set $\la_{kk} :=1$ and $\la_{jk} := \la_{kj}^{-1}$ for $j< k$.
This gives rise to a multiplicatively skew-symmetric 
matrix $\lab := (\la_{kj}) \in M_N(\kx)$.

The CGL extension $R$ is called {\em{torsionfree}} if the subgroup $\langle \la_{kj} \mid k,j \in [1,N] \rangle$
of $\kx$  is torsionfree.
Define the  \emph{rank} of $R$ by 
\begin{equation}
\label{rk}
\rk(R) := \{ k \in [1,N] \mid \de_k = 0 \} \in \Zset_{ > 0}
\end{equation}
(cf.~\cite[Eq. (4.3)]{GY}).

Denote the character group of the torus $\HH$ by $\xh$ and express this group additively. The action of 
$\HH$ on $R$ gives rise to an $\xh$-grading of $R$. The $\HH$-eigenvectors 
are precisely the nonzero homogeneous elements with respect to this grading. 
We denote the $\HH$-eigenvalue of a nonzero homogeneous element $u \in R$ by $\chi_u$. In other words, $\chi_u = \xh\mbox{-deg}(u)$ in terms of the $\xh$-grading.
\ede 

\bpr{LLRHUFD} {\rm\cite[Proposition 3.2]{LLR}} Every CGL extension is an $\HH$-UFD, with $\HH$ as in the definition.
\epr

The sets of homogeneous prime elements in the subalgebras $R_k$ of a CGL extension \eqref{itOre} were characterized in \cite{GY}. The statement of the result involves the standard predecessor and successor functions, $p = p_\eta$ and $s = s_\eta$, of a function $\eta : [1,N] \to \Zset$, defined as follows:
\begin{equation}
\label{pred.succ}
\begin{aligned}
p(k) &= \max \{ j <k \mid \eta(j) = \eta(k) \}, \\
s(k) &= \min \{ j > k \mid \eta(j) = \eta(k) \}, 
\end{aligned}
\end{equation}
where $\max \varnothing = -\infty$ and $\min \varnothing = +\infty$. Define corresponding order functions $O_\pm : [1,N] \rightarrow \Nset$ by
\begin{equation}
\label{O-+}
\begin{aligned}
O_-(k) &:= \max \{ m \in \Nset \mid p^m(k) \ne -\infty \},  \\
O_+(k) &:= \max \{ m \in \Nset \mid s^m(k) \ne +\infty \}.
\end{aligned}
\end{equation}

\bth{CGL} {\rm\cite[Theorem 4.3, Corollary 4.8]{GY}} Let $R$ be a CGL extension of length $N$ as 
in \eqref{itOre}. There exist a function $\eta : [1,N] \to \Zset$ 
whose range has cardinality $\rk(R)$ and elements
$$
c_k \in R_{k-1}, \; \; \forall k \in [2,N] \; \; 
\mbox{with} \; \; p(k) \neq - \infty
$$
such that the elements $y_1, \ldots, y_N \in R$, recursively defined by 
\begin{equation}
\label{y}
y_k := 
\begin{cases}
y_{p(k)} x_k - c_k, &\mbox{if} \; \;  p(k) \neq - \infty \\
x_k, & \mbox{if} \; \; p(k) = - \infty,  
\end{cases}
\end{equation}
are homogeneous and have the property that for every $k \in [1,N]$,
\begin{equation}
\label{prime-elem}
\{y_j \mid j \in [1,k] , \, s(j) > k \}
\end{equation}
is a list of the homogeneous prime elements of $R_k$ up to scalar multiples.

The elements $y_1, \ldots, y_N \in R$ with these properties are unique.
The function $\eta$ satisfying the above 
conditions is not unique, but the partition of $[1,N]$ into a disjoint 
union of the level sets of $\eta$ is uniquely determined by $R$, as are the predecessor and successor functions $p$ and $s$.
The function $p$ has the property that $p(k) = - \infty$
if and only if $\de_k =0$.

Furthermore, the elements $y_k$ 
of $R$ satisfy 
\begin{equation}
\label{prime-comm}
y_k x_j = \al_{jk}^{-1} x_j y_k, \; \; 
\forall k,j \in [1,N], \; s(k) = +\infty,
\end{equation}
where 
$$
\alpha_{jk} := \prod_{m=0}^{O_-(k)} \la_{j,p^m(k)} , \; \; \forall j,k \in [1,N].
$$
\eth

The uniqueness of the level sets of $\eta$ was not stated in \cite[Theorem 4.3]{GY}, 
but it follows at once from \cite[Theorem 4.2]{GY}. This uniqueness immediately implies the uniqueness of $p$ and $s$. 
In the setting of the theorem, the rank of $R$ is also given by
\begin{equation}
\label{rankR}
\rk(R) = |\{ j \in [1,N] \mid s(j) = +\infty \}|
\end{equation}
\cite[Eq. (4.3)]{GY}.

\bpr{prime-x} Let $R$ be a CGL extension of length $N$ as 
in \eqref{itOre}. The following are equivalent 
for an integer $i \in [1,N]$: 

{\rm (a)} The integer $i$ satisfies $\eta^{-1}(\eta(i)) = \{i \}$ for the function $\eta$ from Theorem {\rm\ref{tCGL}}.

{\rm (b)} The element $x_i$ is prime in $R$.

{\rm (c)} The element $x_i$ satisfies $x_i x_j = \la_{ij} x_j x_i$ for all $j \in [1,N]$.
\epr

We will denote by $P_x(R)$ the set of integers $i \in [1,N]$ satisfying the 
conditions (a)--(c).

\begin{proof}
Denote by $A$, $B$, and $C$ the sets of integers $i$ occurring in parts (a), (b) and (c). 
We will prove that $A \subseteq B \subseteq C \subseteq A$. 
The inclusion $A \subseteq B \cap C$ follows at once from \thref{CGL} and Eq. \eqref{prime-comm}. 
Moreover, $B \subseteq A$ because of \eqref{prime-comm}, since if $x_i$ is prime it must be a scalar multiple of $y_i$.

Let $i \in C$. Then $x_i x_j = \sig_i(x_j) x_i$ for all $j < i$, whence $\de_i =0$ and so $p(i) = -\infty$. Thus, $\eta^{-1}(\eta(i)) \subseteq [i,N]$. 
Assume that $\eta^{-1}(\eta(i)) \neq \{ i \}$, which implies $s(i) \neq + \infty$. Set $j:=s(i) \in \eta^{-1}(\eta(i))$. Then $p(j) = i \ne -\infty$ and so $\de_j(x_i) \neq 0$ by \cite[Proposition 4.7 (b)]{GY}, which contradicts the equality $x_j x_i = \la_{ij}^{-1} x_i x_j = \sig_j(x_i) x_j$. 
Therefore $\eta^{-1}(\eta(i)) = \{ i \}$ and $i \in A$.     
\end{proof}

One can show that the conditions in \prref{prime-x} (a)--(c) are equivalent to $x_i$ being a normal 
element of $R$. 

Recall that \emph{quantum tori} and \emph{quantum affine space algebras} over $\KK$ are defined by
\begin{align*}
\TT_\pb &= \OO_\pb((\kx)^N) := \KK \lcor Y_1^{\pm 1}, \ldots, Y_N^{\pm 1} \mid Y_k Y_j = p_{kj} Y_j Y_k, \; \forall k,j \in [1,N] \rcor , \\
\AA_\pb &= \OO_\pb(\KK^N) := \KK \lcor Y_1, \ldots, Y_N \mid Y_k Y_j = p_{kj} Y_j Y_k, \; \forall k,j \in [1,N] \rcor,
\end{align*}
for any multiplicatively skew-symmetric matrix $\pb = (p_{ij}) \in M_N(\kx)$.

\bpr{AinRinT}  {\rm\cite[Theorem 4.6]{GY}}
 For any CGL extension $R$ of length $N$, the elements 
$y_1, \dots, y_N$ generate a quantum affine space algebra $\AA$ 
inside $R$. The corresponding quantum torus $\TT$ is naturally embedded 
in $\Fract(R)$ and we have the inclusions
$$  
\AA \subseteq R \subset \TT.
$$
\epr

The algebras $\AA$ and $\TT$ in \prref{AinRinT} are isomorphic to $\AA_\qb$ and $\TT_\qb$, respectively, where by \cite[Eq. (4.17)]{GY} the entries of the matrix $\qb = (q_{ij})$ are given by
\begin{equation}
\label{defq}
q_{kj} = \prod_{m=0}^{O_-(k)} \prod_{l=0}^{O_-(j)} \la_{p^m(k), p^l(j)}, \; \; \forall k,j \in [1,N].
\end{equation}

\bde{defsaturated}
Let $\pb = (p_{ij}) \in M_N(\kx)$ be a multiplicatively skew-symmetric matrix. Define the (skew-symmetric) multiplicative bicharacter $\Om_\pb : \Zset^N \times \Zset^N \rightarrow \kx$ by
$$
\Om_\pb(e_i, e_j) = p_{ij}, \; \; \forall i,j \in [1,N],
$$
where $e_1, \dots, e_N$ denotes the standard basis of $\Zset^N$. The \emph{radical} of $\Om_\pb$ is the subgroup
$$
\rad \Om_\pb := \{ f \in \Zset^N \mid \Om_\pb(f,g) = 1, \; \forall g \in \Zset^N \}
$$
of $\Zset^N$. We say that the bicharacter $\Om_\pb$ is \emph{saturated} if $\Zset^N / \rad \Om_\pb$ is torsionfree, i.e.,
$$
nf \in \rad \Om_\pb \implies f \in \rad \Om_\pb, \; \; \forall n \in \Zset_{>0}, \; f \in \Zset^N .
$$
Carrying the terminology forward, we say that the quantum torus $\TT_\pb$ is \emph{saturated} provided $\Om_\pb$ is saturated.

Finally, we apply this terminology to a CGL extension $R$ via its associated matrix $\lab$, and say that $R$ is \emph{saturated} provided the bicharacter $\Om_{\lab}$ is saturated.
\ede

For example, any torsionfree CGL extension $R$ is saturated, because all values of $\Om_{\lab}$ lie in the torsionfree group $\langle \la_{kj} \mid k,j \in [1,N] \rangle$ in that case.

\ble{RsatT}
Let $R$ be a CGL extension of length $N$ as in \eqref{itOre}, and let $\TT$ be the quantum torus  in Proposition {\rm\ref{pAinRinT}}. Then $R$ is saturated if and only if $\TT$ is saturated.
\ele

\begin{proof} In view of \eqref{defq}, $\Om_\qb(e_k, e_j) = q_{kj} = \Om_{\lab}(\ol{e}_k, \ol{e}_j)$ for all $k,j \in [1,N]$, where
$$
\ol{e}_i := e_i + e_{p(i)} + \cdots + e_{O_-(i)}, \; \; \forall i \in [1,N].
$$
Since $\ol{e}_1, \dots, \ol{e}_N$ is a basis for $\Zset^N$, it follows that $\Om_{\lab}$ is saturated if and only if $\Om_\qb$ is saturated.
\end{proof}

Continue to let $R$ be a CGL extension of length $N$ as in \eqref{itOre}.
Denote by $\NN(R)$ the unital subalgebra of $R$ generated by its homogeneous 
prime elements $y_k$, $k \in [1,N]$, $s(k) = + \infty$. By \cite[Proposition 2.6]{GY},
$\NN(R)$ coincides with the unital subalgebra of $R$ generated by all normal elements of $R$. As in \leref{IcapER}, denote by $E(R)$ the multiplicative subset of $R$ generated by the homogeneous prime elements of $R$. In the present situation, $E(R)$ is also generated by the set $\kx \sqcup \{ y_k \mid k \in [1,N], \; s(k) = + \infty \}$. It is an Ore set in $R$ and $\NN(R)$ since it is generated by elements which are normal in both algebras. Note that $\NN(R)[E(R)^{-1}] \subseteq R[E(R)^{-1}] \subseteq \TT$, where $\TT$ is the torus of \prref{AinRinT}.

\bpr{center} The center of the quantum torus $\TT$ in Proposition {\rm\ref{pAinRinT}} coincides with the center of $R[E(R)^{-1}]$ and is contained in $\NN(R)[E(R)^{-1}]$, i.e.,
\begin{equation}
\label{ZTinNE-1}
Z(\TT) = Z \bigl( R[E(R)^{-1}] \bigr) =  \{ z \in \NN(R)[E(R)^{-1}] \mid zx = xz, \; \forall x \in R \}.
\end{equation}
\epr

\begin{proof} It is clear that $Z \bigl( R[E(R)^{-1}] \bigr) \subseteq Z(\TT)$, because these centers consist of the elements in $R[E(R)^{-1}]$ and $\TT$ that commute with all elements of $R$, and that the set on the right hand side of \eqref{ZTinNE-1} is contained in $Z \bigl( R[E(R)^{-1}] \bigr)$. Hence, it suffices to show that $Z(\TT) \subseteq \NN(R)[E(R)^{-1}]$. 

Recall that the center of any quantum torus equals the linear span of the central Laurent monomials in its generators. If $m$ is a central Laurent monomial in the generators $y_1^{\pm1}, \dots, y_N^{\pm1}$ of $\TT$, then $m$ is an $\HH$-eigenvector and 
$$
I := \{ r \in R \mid mr \in R \}
$$
is a nonzero $\HH$-stable ideal of $R$. By \leref{IcapER}, there exists $c \in I \cap E(R)$, and $m = ac^{-1}$ for some $a \in R$. Since $m$ centralizes $R$ and $c$ normalizes it, the element $a = mc$ is normal in $R$. Hence, $a \in \NN(R)$, and we conclude that $m \in \NN(R)[E(R)^{-1}]$. There\-fore $Z(\TT) \subseteq \NN(R)[E(R)^{-1}]$, as required.
\end{proof}

\bde{daut}
We will say that an automorphism $\psi$ of a CGL extension $R$ as in \eqref{itOre} is \emph{diagonal} provided $x_1,\dots,x_N$ are eigenvectors for $\psi$. Set
$$
\DAut(R) := \{ \text{diagonal automorphisms of} \; R \},
$$
a subgroup of $\Aut(R)$. 

In particular, the group $\{ (h\cdot) \mid h \in \HH \}$ is contained in $\DAut(R)$. It was shown in \cite[\S 5.2 and Theorem 5.2]{GY} that $\DAut(R)$ is naturally isomorphic to a $\KK$-torus of rank equal to $\rank(R)$, exhibited as a closed connected subgroup of the torus $(\kx)^N$. This allows us to think of $\DAut(R)$ as a torus, and to replace $\HH$ by $\DAut(R)$ if desired. A description of this torus, as a specific subgroup of $(\kx)^N$, was established in \cite[Theorem 5.5]{GY}. Finally, it follows from \cite[Corollary 5.4]{GY}  that for any nonzero normal element $u \in R$, there exists $\psi \in \DAut(R)$ such that $ua = \psi(a)u$ for all $a \in R$.
\ede

\subsection{Symmetric CGL extensions} \label{CGL.symm}

For a CGL extension $R$ as in \eqref{itOre} and $j,k \in [1,N]$, denote by $R_{[j,k]}$ the unital subalgebra of
$R$ generated by $\{ x_i \mid j \le i \le k \}$. So, $R_{[j,k]} = \KK$ if $j \nleq k$. 

\bde{symmetric} We call a CGL extension $R$ of length $N$ as in 
\deref{CGL} {\em{symmetric}} if the following two conditions hold:
\begin{enumerate}
\item[(i)] For all $1 \leq j < k \leq N$,
$$
\de_k(x_j) \in R_{[j+1, k-1]}.
$$
\item[(ii)] For all $j \in [1,N]$, there exists $h^*_j \in \HH$ 
such that 
$$
h^*_j \cdot x_k = \la_{kj}^{-1} x_k = \la_{jk} x_k, \; \; \forall 
k \in [j+1, N]
$$
and $h^*_j \cdot x_j = \la^*_j x_j$ for some $\la^*_j \in \kx$ which is not a root of unity.
\end{enumerate}
\ede

For example, all quantum Schubert cell algebras $U^+[w]$ are symmetric CGL extensions, 
cf.~\exref{qsc} below. 

Given a symmetric CGL extension $R$ as in \deref{symmetric}, set 
$$
\sig^*_j := (h^*_j \cdot) \in \Aut (R), \; \; \forall j \in [1,N-1].
$$
Then for all $j \in [1,N-1]$, 
the inner $\sig^*_j$-derivation on $R$ 
given by $a \mapsto x_j a - \sig^*_j(a) x_j$ restricts to a $\sig^*_j$-derivation 
$\de^*_j$ of $R_{[j+1, N]}$. It is given by
$$
\de^*_j(x_k): = x_j x_k - \la_{jk} x_k x_j = - \la_{jk} \de_k(x_j), \; \; 
\forall k \in [j+1, N].
$$ 
For all $1 \leq j < k \leq N$, $\sig_k$ and $\de_k$ preserve $R_{[j,k-1]}$
and $\sig^*_j$ and $\de^*_j$ preserve $R_{[j+1,k]}$. This gives 
rise to the skew polynomial extensions
\begin{equation}
\label{step1}
R_{[j,k]} = R_{[j,k-1]}[x_k; \sig_k, \de_k] \; \; 
\mbox{and} \; \; 
R_{[j,k]} = R_{[j+1,k]} [ x_j; \sig^*_j, \de^*_j].
\end{equation}
In particular, it follows that $R$ has an iterated Ore extension presentation with the variables $x_k$ in descending order:
\begin{equation}
\label{CGLdescending}
R = \KK[x_N] [x_{N-1}; \sig^*_{N-1}, \de^*_{N-1}] \cdots [x_1; \sig^*_1, \de^*_1].
\end{equation}
This is the reason for the name ``symmetric".

Denote the following subset of the symmetric group $S_N$:
\begin{multline}
\label{tau}
\Xi_N := \{ \tau \in S_N \mid 
\tau(k) = \max \, \tau( [1,k-1]) +1 \; \;
\mbox{or} 
\\
\tau(k) = \min \, \tau( [1,k-1]) - 1, 
\; \; \forall k \in [2,N] \}.
\end{multline}
In other words, $\Xi_N$ consists of those $\tau \in S_N$ 
such that $\tau([1,k])$ is an interval for all $k \in [2,N]$. 
For each $\tau \in \Xi_N$, we have the iterated 
Ore extension presentation
\begin{equation}
\label{tauOre}
R = \KK [x_{\tau(1)}] [x_{\tau(2)}; \sig''_{\tau(2)}, \de''_{\tau(2)}] 
\cdots [x_{\tau(N)}; \sig''_{\tau(N)}, \de''_{\tau(N)}],
\end{equation}
where $\sig''_{\tau(k)} := \sig_{\tau(k)}$ and 
$\de''_{\tau(k)} := \de_{\tau(k)}$ if 
$\tau(k) = \max \, \tau( [1,k-1]) +1$, while 
$\sig''_{\tau(k)} := \sig^*_{\tau(k)}$ and 
$\de''_{\tau(k)} := \de^*_{\tau(k)}$ if 
$\tau(k) = \min \, \tau( [1,k-1]) -1$.

\bpr{tauOre} 
\cite[Remark 6.5]{GY} 
For every symmetric CGL extension $R$ of length $N$ and any $\tau \in \Xi_N$,
the iterated Ore extension presentation \eqref{tauOre} of $R$ 
is a CGL extension presentation for the 
same choice of $\KK$-torus $\HH$, and the associated elements 
$h''_{\tau(1)}, \ldots, h''_{\tau(N)} \in \HH$ required by Definition {\rm\ref{dCGL}(iii)} 
are given by $h''_{\tau(k)} = h_{\tau(k)}$ if $\tau(k) = \max \, \tau( [1,k-1]) +1$  
and $h''_{\tau(k)} = h^*_{\tau(k)}$ if $\tau(k) = \min \, \tau( [1,k-1]) -1$.
\epr

It follows from \prref{tauOre} that in the given situation, $\sig''_{\tau(k)}(x_{\tau(j)}) = \la_{\tau(k),\tau(j)} x_{\tau(j)}$ for $1\le j < k \le N$. Hence, the $\lab$-matrix for the presentation \eqref{tauOre} is the matrix 
\begin{equation}
\label{labtau}
\lab_\tau := (\la_{\tau(k), \tau(j)}).
\end{equation}

If $R$ is a symmetric CGL extension of length $N$ and $\tau \in \Xi_N$, we write $y_{\tau, 1}, \ldots, y_{\tau, N}$ for the $y$-elements obtained from applying \thref{CGL} to the CGL extension presentation \eqref{tauOre}. \prref{AinRinT} then shows that $y_{\tau, 1}, \ldots, y_{\tau, N}$ generate a quantum affine space algebra $\AA_\tau$ 
inside $R$, the corresponding quantum torus $\TT_\tau$ is naturally embedded 
in $\Fract(R)$, and we have the inclusions
$$ 
\AA_\tau \subseteq R \subset \TT_\tau.
$$

\bpr{sat} If $R$ is a saturated symmetric CGL extension of length $N$, then the quantum tori $\TT_\tau$ are saturated, for all $\tau \in \Xi_N$.
\epr

\begin{proof} Let $\tau \in \Xi_N$, and recall \eqref{labtau}. It follows that
$$
\Om_{\lab_\tau}(f,g) = \Om_{\lab}( \tau \cdot f, \tau \cdot g ), \; \; \forall f,g \in \Zset^N,
$$ 
where we identify $\tau$ with the corresponding permutation matrix in $GL_N(\Zset)$ and write elements of $\Zset^N$ as column vectors. Since $\Om_{\lab}$ is saturated by hypothesis, it follows immediately that $\Om_{\lab_\tau}$ is saturated. Applying \leref{RsatT} to the presentation \eqref{tauOre}, we conclude that $\TT_\tau$ is saturated.
\end{proof}

\sectionnew{Nakayama automorphisms of iterated Ore extensions}
\label{Nakayama}

Every iterated Ore extension $R$ over $\KK$ is a twisted Calabi-Yau algebra (see \deref{tCY} and \coref{CGLCY}), and as such has an associated Nakayama automorphism, which is unique in this case because the inner automorphisms of $R$ are trivial. Our main aim is to determine this automorphism $\nu$ when $R$ is a symmetric CGL extension. In that case, we show that $\nu$ is the restriction to $R$ of an inner automorphism $u^{-1}(-)u$ of $\Fract(R)$, where $u = u_1 \cdots u_n$ for a list $u_1,\dots, u_n$ of the homogeneous prime elements of $R$ up to scalar multiples. On the way, we formalize a technique of Liu, Wang, and Wu \cite{LWW} and use it to give a formula for $\nu$ in a more general symmetric situation, where we show that each standard generator of $R$ is an eigenvector for $\nu$ and determine the eigenvalues.

Recall that the \emph{right twist} of a bimodule $M$ over a ring $R$ by an automorphism $\nu$ of $R$ is the $R$-bimodule $M^\nu$ based on the left $R$-module $M$ and with right $R$-module multiplication $*$ given by $m*r = m \nu(r)$ for $m \in M$, $r \in R$.

\bde{tCY} A $\KK$-algebra $R$ is \emph{$\nu$-twisted Calabi-Yau of dimension $d$}, where $\nu$ is an automorphism of $R$ and $d \in \Znn$, provided
\begin{enumerate}
\item[(i)] $R$ is \emph{homologically smooth}, meaning that as a module over $R^e := R \otimes_\KK R^{\text{op}}$, it has a finitely generated projective resolution of finite length.
\item[(ii)] $\Ext^i_{R^e}(R,R^e) \cong \begin{cases} 0  &(\text{if\;} i\ne d)  \\  R^\nu  &(\text{if\;} i=d)  \end{cases}$ as $R^e$-modules.
\end{enumerate}
When these conditions hold, $\nu$ is called the \emph{Nakayama automorphism} of $R$. It is unique up to an inner automorphism. The case of a Calabi-Yau algebra in the sense of Ginzburg \cite{Gin} is recovered when $\nu$ is inner.
\ede

\bth{LWWthm}  {\rm(Liu-Wang-Wu \cite[Theorem 3.3]{LWW})}
Let $B$ be a $\nu_0$-twisted Calabi-Yau algebra of dimension $d$, and let $R := B[x;\sig,\de]$ be an Ore extension of $B$. Then $R$ is a $\nu$-twisted Calabi-Yau algebra of dimension $d+1$, where $\nu$ satisfies the following conditions:
\begin{enumerate}
\item[\rm(a)] $\nu|_B = \sig^{-1}\nu_0$.
\item[\rm(b)] $\nu(x)= ux+b$ for some unit $u\in B$ and some $b\in B$.
\end{enumerate}
\eth

\bco{CGLCY}
Every iterated Ore extension of length $N$ over $\KK$ is a twisted Calabi-Yau algebra of dimension $N$.
\eco

Note that the only units in an iterated Ore extension $R$ over $\KK$ are scalars, and so the only inner automorphism of $R$ is the identity. Hence, the Nakayama automorphism of $R$ is unique.

Liu, Wang, and Wu gave several examples in \cite{LWW} for which the Nakayama automorphism can be completely pinned down by \thref{LWWthm}. These examples are iterated Ore extensions which can  be rewritten as iterated Ore extensions with the original variables in reverse order. We present a general result of this form in the following subsection, and apply it to symmetric CGL extensions in \S \ref{sCGL}.

\subsection{Nakayama automorphisms of reversible iterated Ore extensions}
\label{revOre}

\bde{revitOre}
Let $R$ be an iterated Ore extension of length $N$ as in \eqref{itOre}. We shall say that $R$ (or, more precisely, the presentation \eqref{itOre}) is \emph{diagonalized} if there are scalars $\la_{kj} \in \kx$ such that $\sig_k(x_j) = \la_{kj} x_j$ for all $1 \leq j < k \leq N$. When $R$ is diagonalized, we extend the $\la_{kj}$ to a multiplicatively skew-symmetric matrix just as in the CGL case.

A diagonalized iterated Ore extension $R$ is called \emph{reversible} provided
there is a second iterated Ore extension presentation
\begin{equation}  \label{itOre2}
R = \KK [x_N] [x_{N-1}; \sig^*_{N-1}, \de^*_{N-1}] \cdots [x_1; \sig^*_1, \de^*_1]
\end{equation}
such that $\sig^*_j(x_k) = \la_{jk} x_k$ for all $1 \le j < k \le N$. 
\ede

Every symmetric CGL extension is a reversible diagonalized iterated Ore extension, by virtue of the presentation \eqref{CGLdescending}.

For any iterated Ore extension $R$ as in \eqref{itOre}, we define the subalgebras $R_{[j,k]}$ of $R$ just as in \S \ref{CGL.symm}.

\ble{revcond}
Let $R$ be a diagonalized iterated Ore extension of length $N$ as in Definition {\rm\ref{drevitOre}}. Then $R$ is reversible if and only if
\begin{equation}
\label{revcond}
\de_k(x_j) \in R_{[j+1,k-1]}, \; \; \forall 1\le j<k\le N .
\end{equation}
\ele

\begin{proof} 
Assume first that $R$ is reversible, and let  $1\le j<k\le N$. From the structure of the iterated Ore extensions \eqref{itOre} and \eqref{itOre2}, we see that
$$\de_k(x_j) \in R_{[1,k-1]} \qquad\quad \text{and} \qquad\quad \de^*_j(x_k) \in R_{[j+1,N]}.$$
Since $R$ is diagonalized, we also have
$$\de^*_j(x_k) = x_j x_k - \la_{jk} x_k x_j = - \la_{jk} \bigl( x_k x_j - \sig_k(x_j) x_k \bigr) = - \la_{jk} \de_k(x_j),$$
and thus $\de_k(x_j) \in R_{[j+1,N]}$. Since $R_{[1,k-1]}$ and $R_{[j+1,N]}$ are iterated Ore extensions with PBW bases $\{ x_1^\bullet \cdots x_{k-1}^\bullet \}$ and $\{ x_{j+1}^\bullet \cdots x_N^\bullet \}$, respectively, it follows that $\de_k(x_j) \in R_{[j+1,k-1]}$, verifying \eqref{revcond}.

Conversely, assume that \eqref{revcond} holds. We establish the following by a  downward induction on $l \in [1,N]$:
\begin{enumerate}
\item[(a)] The monomials
\begin{equation}  \label{monom.lN}
x_l^{a_l} \cdots x_N^{a_N}, \; \; \forall a_l,\dots,a_N \in \Znn
\end{equation}
form a basis of $R_{[l,N]}$.
\item[(b)] $R_{[l,N]} = R_{[l+1,N]} [x_l; \sig^*_l, \de^*_l]$
for some automorphism $\sig^*_l$ and $\sig^*_l$-derivation $\de^*_l$ of $R_{[l+1,N]}$, such that $\sig^*_l(x_k) = \la_{lk} x_k$ for all $k \in [l+1,N]$.
\end{enumerate}
When $l=N$, both (a) and (b) are clear, since $R_{[N,N]} = \KK[x_N]$ and $R_{[N+1,N]} = \KK$.

Now let $l \in [1,N-1]$ and assume that (a) and (b) hold for $R_{[l+1,N}]$. For $k \in [l+1,N]$, it follows from \eqref{revcond} that 
$$x_k x_l - \la_{kl} x_l x_k = \de_k(x_l) \in R_{[l+1,k-1]} \subset R_{[l+1,N]}.$$
Consequently, we see that
\begin{equation}  \label{innerl}
R_{[l+1,N]} + x_l R_{[l+1,N]} = R_{[l+1,N]} + R_{[l+1,N]} x_l \,.
\end{equation}
In particular, \eqref{innerl} implies that $\sum_{a=0}^\infty x_l^a R_{[l+1,N]}$ is a subalgebra of $R$. In view of our induction hypothesis, it follows that the monomials \eqref{monom.lN} span $R_{[l,N]}$. Consequently, they form a basis, since they are part of the standard PBW basis for $R$. This establishes (a) for $R_{[l,N]}$. 

Given the above bases for $R_{[l,N]}$ and $R_{[l+1,N]}$, we see that $R_{[l,N]}$ is a free right $R_{[l+1,N]}$-module with basis $(1,x_l,x_l^2,\dots)$. Via \eqref{innerl} and an easy induction on degree, we confirm that $R_{[l,N]}$ is also a free left $R_{[l+1,N]}$-module with the same basis. A final application of \eqref{innerl} then yields $R_{[l,N]} = R_{[l+1,N]} [x_l; \sig^*_l, \de^*_l]$
for some automorphism $\sig^*_l$ and $\sig^*_l$-derivation $\de^*_l$ of $R_{[l+1,N]}$. For $k \in [l+1,N]$, we have
$$x_l x_k - \la_{lk} x_k x_l = - \la_{lk} \bigl( x_k x_l - \sig_k(x_l) x_k \bigr) = - \la_{lk} \de_k(x_l) \in R_{[l+1,k-1]} \subset R_{[l+1,N]} \,,$$
from which it follows that $\sig^*_l(x_k) = \la_{lk} x_k$. Thus, (b) holds for $R_{[l,N]}$.

Therefore, the induction works. Combining statements (b) for $l=N, \dots, 1$, we conclude that $R$ is reversible.
\end{proof}

\bth{Nak.rditOre}
Let $R = \KK[x_1][x_2; \sig_2, \delta_2] \cdots [x_N; \sig_N, \delta_N]$ be a reversible, diagonalized iterated Ore extension over $\KK$, let $\nu$ be the Nakayama automorphism of $R$, and let $(\la_{jk}) \in M_N(\kx)$ be the multiplicatively antisymmetric matrix such that
$\sig_k(x_j) = \la_{kj} x_j$ for all $1 \leq j < k \leq N$. Then
\begin{equation}  \label{Nak.form}
\nu(x_k) = \biggl( \prod_{j=1}^N \la_{kj} \biggr) x_k, \; \; \forall k \in [1,N].
\end{equation}
\eth

\begin{proof} In case $N=1$, the algebra $R$ is a polynomial ring $\KK[x_1]$. Then $R$ is Calabi-Yau (e.g., as in \cite[Example 13]{Far}), i.e., $\nu$ is the identity. Thus, the theorem holds in this case.

Now let $N\ge2$, and assume the theorem holds for all reversible, diagonalized iterated Ore extensions of length less than $N$. It is clear from the original and the reversed iterated Ore extension presentations of $R$ that $R_{N-1}$ and $R_{[2,N]}$ are diagonalized iterated Ore extensions, and it follows from \leref{revcond} that $R_{N-1}$ and $R_{[2,N]}$ are reversible.
If $\nu_0$ denotes the Nakayama automorphism of $R_{N-1}$, then the inductive statement together with \thref{LWWthm} gives us
\begin{equation}  \label{Nak.1N-1}
\nu(x_k) = \sig_N^{-1} \nu_0(x_k) = \la_{Nk}^{-1} \biggl( \prod_{j=1}^{N-1} \la_{kj} \biggr) x_k = \biggl( \prod_{j=1}^N \la_{kj} \biggr) x_k, \; \; \forall k \in [1,N-1].
\end{equation}
Similarly, if $\nu_1$ denotes the Nakayama automorphism of $R_{[2,N]}$, we obtain
\begin{equation}  \label{Nak.2N}
\nu(x_k) = (\sig^*_1)^{-1} \nu_1(x_k) = \la_{1k}^{-1} \biggl( \prod_{j=2}^N \la_{kj}  \biggr) x_k = \biggl( \prod_{j=1}^N \la_{kj} \biggr) x_k, \; \; \forall k \in [2,N].
\end{equation}
The formulas \eqref{Nak.1N-1} and \eqref{Nak.2N} together yield \eqref{Nak.form}, establishing the induction step.
\end{proof}

In particular, \thref{Nak.rditOre} immediately determines the Nakayama automorphisms of the multiparameter quantum affine spaces $\OO_{\qb}(\KK^N)$, as in \cite[Proposition 4.1]{LWW}, and those of the Weyl algebras $A_n(\KK)$ \cite[Remark 4.2]{LWW}.

\bex{OqMtn}
Let $R := \OO_q(M_{t,n}(\KK))$ be the standard quantized coordinate ring of $t\times n$ matrices over $\KK$, with $q \in \kx$, generators $X_{ij}$ for $i\in [1,t]$, $j\in [1,n]$, and relations
\begin{align*}
X_{ij}X_{im} &= qX_{im}X_{ij}  &X_{ij}X_{lj} &= qX_{lj}X_{ij} \\
X_{im}X_{lj} &= X_{lj}X_{im}  &X_{ij}X_{lm}-X_{lm}X_{ij} &= (q-q^{-1})X_{im}X_{lj}
\end{align*}
for $i<l$ and $j<m$. (We place no restriction on $q$ for now.) It is well known that $R$ has an iterated Ore extension presentation with variables $X_{ij}$ listed in lexicographic order, that is,
\begin{equation}
\label{OqMtnCGL}
\begin{gathered}
R := \KK[x_1][x_2; \sig_2, \delta_2] \cdots [x_N; \sig_N, \delta_N], \; \; N := tn,  \\
x_{(i-1)n+j} := X_{ij}, \; \; \forall i\in [1,t], \; j\in [1,n].
\end{gathered}
\end{equation}
It is clear that $R$ is diagonalized. Since $R$ also has an iterated Ore extension presentation with the  $X_{ij}$ in reverse lexicographic order, one easily checks that $R$ is thus reversible. 

The scalars $\la_{(i-1)n+j, \, (l-1)n+m}$ are equal to $1$ except in the following cases:
\begin{align*}
\la_{(i-1)n+j, \, (i-1)n+m} &= q^{-1} \quad(m<j)  &\la_{(i-1)n+j, \, (i-1)n+m} &= q \quad(m>j)  \\
\la_{(i-1)n+j, \, (l-1)n+j} &= q^{-1} \quad(l<i)  &\la_{(i-1)n+j, \, (l-1)n+j} &= q \quad(l>i).
\end{align*}
In view of \thref{Nak.rditOre}, we thus find that the Nakayama automorphism $\nu$ of $R$ is given by the rule
$$\nu(X_{ij}) = q^{t+n-2i-2j+2} X_{ij}$$
for $i\in [1,t]$, $j\in [1,n]$.

Let us only consider the multiparameter version of $R$ in the $n\times n$ case. This is the $\KK$-algebra $R' := \OO_{\la,\pb}(M_{n}(\KK))$, where $\la \in \KK\backslash \{0,1\}$ and $\pb$ is a multiplicatively skew-symmetric $n\times n$ matrix over $\kx$, with generators $X_{ij}$ for $i,j \in [1,n]$ and relations
$$
X_{l m}X_{ij} = \begin{cases} p_{l i}p_{jm}X_{ij}X_{l m} +
(\lambda -1)p_{l i}X_{im}X_{l j}&\quad (l >i,\
m>j)\\ 
\lambda p_{l i}p_{jm}X_{ij}X_{l m}&\quad (l >i,\ m\le j)\\ 
p_{jm}X_{ij}X_{l m}&\quad (l
=i,\ m>j).\end{cases}
$$
Iterated Ore extension presentations of $R'$ are well known, and as above, we see that $R'$ is diagonalized and reversible. It follows from \thref{Nak.rditOre} that
$$
\nu(X_{ij}) = \biggl( \prod_{l=1}^n p_{il}^n \biggr) \biggl( \prod_{m=1}^n p_{mj}^n \biggr) \la^{n(i-j-1)+i+j-1} X_{ij}, \; \; \forall i,j \in [1,n].
$$
\eex

\subsection{Nakayama automorphisms of symmetric CGL extensions}
\label{sCGL}

As noted above, any symmetric CGL extension is reversible and diagonalized, so \thref{Nak.rditOre} provides a formula for its Nakayama automorphism. We prove that in this case, the Nakayama automorphism arises from the action of a normal element, as follows.

\bth{Nak.symCGL2}
Let $R$ be a symmetric CGL extension of length $N$ as in Definition {\rm\ref{dsymmetric}} and $\nu$ its Nakayama automorphism. Let $u_1,\dots,u_n$ be a complete, irredundant list of the homogeneous prime elements of $R$ up to scalar multiples, and set $u = u_1\cdots u_n$. Then $\nu$ satisfies {\rm(}and is determined by\/{\rm)} the following condition:
\begin{equation}  \label{uxnu}
au = u\nu(a), \; \; \forall a \in R.
\end{equation}
\eth

\begin{proof} Replacing the $u_i$ by scalar multiples of these elements has no effect on \eqref{uxnu}. Thus, we may assume that, in the notation of \thref{CGL}, 
$$
\{u_1,\dots,u_n\} = \{ y_l \mid l\in [1,N],\; s(l) = +\infty \} .
$$
Hence, \eqref{prime-comm} implies that
\begin{equation}  \label{xk.u}
x_k u = \be_k u x_k \;\; \text{with} \;\; \be_k := \prod_{l \in [1,N],\; s(l) = +\infty} \alpha_{kl}, \; \; \forall k \in [1,N].
\end{equation} 
As $l$ runs through the elements of $[1,N]$ with $s(l) = +\infty$ and $m$ runs from $0$ to $O_-(l)$, the numbers $p^m(l)$ run through the elements of $[1,N]$ exactly once each. Hence,
\begin{equation}  \label{expandbeta}
\beta_k =  \prod_{l \in [1,N],\; s(l) = +\infty} \prod_{m=0}^{O_-(l)} \la_{k,p^m(l)} = \prod_{j=1}^N \la_{kj} \;.
\end{equation}
In view of \thref{Nak.rditOre}, we obtain from \eqref{xk.u} and \eqref{expandbeta} that $x_k u = u \nu(x_k)$ for all $k \in [1,N]$. The relation \eqref{uxnu} follows.
\end{proof}

\bex{OqMtn2} Return to $R := \Oq(M_{t,n}(\KK))$ as in \exref{OqMtn}, and assume that $q$ is not a root of unity. There is a rational action of the torus $\HH := (\kx)^{t+n}$ on $R$ by $\KK$-algebra automorphisms such that
$$
(\al_1, \dots, \al_{t+n}) \cdot X_{ij} = \al_i \al_{t+j} X_{ij}, \; \; \forall (\al_1, \dots, \al_{t+n}) \in \HH, \; i \in [1,t], \; j \in [1,n],
$$
and it is well known that $R$ equipped with this action is a CGL extension. It is easily checked that $R$ is symmetric. The function $\eta$ from \thref{CGL} can be chosen so that
$$
\eta((i-1)n+j) = j-i, \; \; \forall i \in [1,t], \; j \in [1,n].
$$ 
The element $y_{(i-1)n+j}$ is the largest solid quantum minor with lower right corner in row $i$, column $j$, that is,
$$
y_{(i-1)n+j} = \bigl[ [i- \min(i,j) +1, i] \bigm| [j- \min(i,j) +1, j] \bigr], \; \; \forall i \in [1,t], \; j \in [1,n],
$$ 
and the list of homogeneous prime elements of $R$, up to scalar multiples, given by \thref{CGL} is
$$
y_n, y_{2n}, \dots, y_{(t-1)n}, y_{(t-1)n+1}, y_{(t-1)n+2}, \dots, y_{tn}.
$$
The product of these $t+n-1$ quantum minors gives the element $u$ that determines the Nakayama automorphism of $R$ as in \thref{Nak.symCGL2}.
\eex

\bex{qsc} Let $\g$ be a simple Lie algebra with set of simple 
roots $\Pi$, Weyl group $W$, and root lattice $Q$, and set $Q^+ := \Zset_{\geq 0} \Pi$. 
For each $\al \in \Pi$, denote by $s_\al \in W$ 
and $\vpi_\al$ the corresponding reflection and 
fundamental weight. 
Denote by $\lcor ., .\rcor$ the $W$-invariant, 
symmetric, bilinear form on $\Qset \Pi$, 
normalized by $\lcor \al, \al \rcor = 2$ for 
short roots $\al$. 
Let $\UU_q(\g)$ be the quantized 
universal enveloping algebra 
of $\g$ over an arbitrary base field 
$\KK$ for a deformation parameter $q \in \KK^*$ 
which is not a root of unity. We will use the 
notation of \cite{Ja}. In particular, we will denote 
the standard generators of $\UU_q(\g)$ by 
$E_\al$, $K_\al^{\pm 1}$, $F_\al$, $\al \in \Pi$.
The subalgebra of $\UU_q(\g)$ generated by 
$\{E_\al \mid \al \in \Pi \}$ will be denoted 
by $\UU_q^+(\g)$. It is naturally $Q^+$-graded 
with $\deg E_\al = \al$ for $\al \in \Pi$.  
For each $w \in W$, De Concini--Kac--Procesi 
and Lusztig defined a graded subalgebra $\UU^+[w]$ 
of $\UU_q^+(\g)$, given by \cite[\S8.21-8.22]{Ja}. 
It is well known that $\UU^+[w]$ is a symmetric CGL extension for the torus 
$\HH := (\KK^*)^{|\Pi|}$ and the action 
$$
t \cdot x := \Big( \prod_{\al \in \Pi} t_\al^{\lcor \al, \ga \rcor} \Big) 
x, \; \; \forall
t= (t_\al)_{\al \in \Pi} \in (\KK^*)^{|\Pi|}, \; x \in \UU_q^+(\g)_\ga, \; \ga \in Q^+.
$$
Here and below, for a
$Q^+$-graded algebra $R$ we denote by 
$R_\ga$ the homogeneous component of $R$ of degree $\ga \in Q^+$.
(Note that the $\HH$-eigenvectors 
in $\UU^+[w]$ are precisely the homogeneous 
elements with respect to the $Q^+$-grading.) 

The algebra $\UU^+[w]$ is a deformation of the universal 
enveloping algebra $\UU(\n_+ \cap w (\n_-))$ 
where $\n_\pm$ are the nilradicals of a pair 
of opposite Borel subalgebras. For each $w \in W$, 
Joseph defined \cite[\S10.3.1]{J} a $Q^+$-graded algebra $S^-_w$
in terms of a localization of the related quantum group 
algebra $R_q[G]$. The grading is given by 
\cite[Eq. (3.22)]{Y-sqg}; here we will omit the trivial 
second component. An explicit $Q^+$-graded 
isomorphism $\vph^-_w : S_w^- \to \UU^+[w]$ 
was constructed in \cite[Theorem 2.6]{Y-sqg}.
Denote the support of $w$
$$
\SS(w) := \{ \al \in \Pi \mid s_\al \leq w \} \subseteq \Pi,
$$
where $\leq$ denotes the Bruhat order on $W$.

For a subset $I \subseteq \Pi$, define the subset 
of dominant integral weights 
$$
P^+_I := \Zset_{\geq 0} \{ \vpi_\al \mid \al \in I \}.
$$
Denote
$$
\rho_I := \sum_{\al \in I} \vpi_\al
$$
which also equals the half-sum of positive roots 
of the standard Levi subalgebra of $\g$ corresponding
to $I$. 

For each $\la \in P^+_{\SS(w)}$, there is a nonzero 
normal element $d^-_{w, \la} \in (S^-_w)_{(1-w)\la}$ given by 
\cite[Eq. (3.29)]{Y-sqg}. It commutes with the elements of 
$S^-_w$ by 
$$
d^-_{w,\la} s = q^{\lcor (w+1)\la, \ga \rcor} s d^-_{w, \la}, \quad
\forall s \in (S^-_w)_\ga, \; \ga \in Q^+.  
$$
We have 
$$
d^-_{w,\la_1} d^-_{w, \la_2} = q^{\lcor \la_1, (1-w) \la_2 \rcor} d^-_{w, \la_1 + \la_2}, 
\quad \forall \la_1, \la_2 \in P^+_{\SS(w)},
$$
cf. \cite[Eq. (3.31)]{Y-sqg}. 
By \cite[Theorem 6.1(ii)]{Y-sqg},  
$$
\{ d^-_{w,\vpi_\al} \mid \al \in \SS(w) \} 
$$
is a list of the homogeneous prime elements of $S^w_-$.
Therefore, up to 
a nonzero scalar multiple the product of the homogeneous prime elements 
of $\UU^+[w]$ is $\vph^-_w(d^-_{w, \rho_{\SS(w)}})$. \thref{Nak.symCGL2} 
implies that the Nakayama automorphism of $\UU^+[w]$ is given by 
$$
\nu(a) = \vph^-_w(d^-_{w, \rho_{\SS(w)}})^{-1} a \vph^-_w(d^-_{w, \rho_{\SS(w)}}), 
\quad \forall a \in \UU^+[w].
$$
Furthermore, the above facts imply that it is also given by 
$$
\nu(a) = q^{-\lcor (w+1) \rho_{\SS(w)} , \ga \rcor} a, 
\quad
\forall a \in \UU^+[w]_\ga, \; \ga \in Q^+.
$$
This is a more explicit form than a previous formula
for the Nakayama automorphism of $\UU^+[w]$ obtained by
Liu and Wu \cite{LW}.
\eex

\sectionnew{Unipotent automorphisms}
\label{Aut}
In this section, we prove a theorem stating that the unipotent 
automorphisms (see \deref{unip}) of a symmetric CGL extension have a very restricted form.
The theorem improves the results in \cite{Y-LL, Y-AD}. 
It is sufficient to classify the full groups of unipotent 
automorphisms of concrete CGL extensions apart from  
examples which have a nontrivial quantum torus factor in a suitable sense. This is illustrated by giving 
a second proof of the Launois--Lenagan 
conjecture \cite{LL1} on automorphisms of square quantum matrix algebras, and by determining the automorphism groups of several other generic quantized coordinate rings.

\subsection{Algebra decompositions of symmetric CGL extensions}
\label{4.1}
Next, we define a unique decomposition of every symmetric CGL extension into 
a crossed product of a symmetric 
CGL extension by a free abelian monoid which has the property that the first term cannot be further so
decomposed. 

Let $R$ be a symmetric CGL extension of length $N$ as in \eqref{itOre}. Recall from 
Section \ref{sectCGL} that $P_x(R) \subseteq[1,N]$ consists of those 
indices $i$ for which $x_i$ is a prime element of $R$. They satisfy 
\begin{equation}
\label{PxR}
x_i x_k = \la_{ik} x_k x_i, \quad \forall k \in [1,N].
\end{equation}
For all $1 \le j < k \le N$, the element
$$
Q_{kj} := x_k x_j - \la_{kj} x_j x_k = \de_k(x_j) \in R_{[j+1, k-1]}
$$
is uniquely a linear combination of monomials $x_{j+1}^{m_{j+1}} \cdots x_{k-1}^{m_{j-1}}$. Of course, $Q_{kj} = 0$ if $k$ or $j$ is in $P_x(R)$.

Denote by $F_x(R)$ the set of those $i \in P_x(R)$ such that $x_i$ does not appear in $Q_{kj}$ (more precisely, no monomial which appears with a nonzero coefficient in $Q_{kj}$ contains a positive power of $x_i$) for any $k, j \in [1,N] \backslash P_x(R)$, $j < k$. 
Let $C_x(R) : = [1,N] \backslash F_x(R)$. The idea for the notation is that 
$F_x(R)$ indexes the set of $x\,$s which will be factored out and $C_x(R)$ 
indexes the set of essential $x\,$s which generate the core of $R$.
Denote the subalgebras 
$$
\CC(R) := \KK \lcor x_k \mid k \in C_x(R) \rcor \qquad \text{and} \qquad
\AA(R) := \KK \lcor x_i \mid i \in F_x(R) \rcor. 
$$
We observe that $R$ is a split extension of either of these subalgebras by a corresponding ideal:
$$
R = \CC(R) \oplus \lcor x_i \mid i \in F_x(R) \rcor = \AA(R) \oplus \lcor x_k \mid k \in C_x(R) \rcor.
$$

Let $C_x(R) = \{ k_1 < k_2 < \ldots < k_t \}$. The algebra $\CC(R)$ is a symmetric
CGL extension of the form
$$
\CC(R) = \KK[x_{k_1}][x_{k_2}; \sig'_{k_2}, \de'_{k_2}] \ldots [x_{k_t}; \sig'_{k_t}, \de'_{k_t}] 
$$ 
where the automorphisms $\sig'_\bullet$, the skew derivations $\de'_\bullet$, and the torus action $\HH$ 
are obtained by restricting those for the CGL extension $R$. The elements $h_\bullet$ and $h^*_\bullet$ 
entering in the definition of a symmetric CGL extension are not changed in going from $R$ to 
$\CC(R)$; we just use a subset of those. The CGL extension $\CC(R)$ will be called the 
\emph{core} of $R$.
The algebra $\AA(R)$ is a quantum affine space algebra with commutation relations 
\begin{equation}
\label{AAR}
x_{i_1} x_{i_2} = \la_{i_1 i_2} x_{i_2} x_{i_1}, \quad \forall i_1, i_2 \in F_x(R). 
\end{equation}
It is a symmetric CGL extension with the restriction of the 
action of $\HH$, but this will not play any role below.

Finally, we can express $R$ as a crossed product
\begin{equation}
\label{decomp}
R = \CC(R) * M,
\end{equation}
where $M$ is a free abelian monoid on $|F_x(R)|$ generators. The actions of these generators on $\CC(R)$ are given by the automorphisms formed  from the commutation relations
\begin{equation}
\label{Rsmash}
x_i x_k = \la_{ik} x_k x_i, \quad \forall i \in F_x(R), \; k \in C_x(R),
\end{equation}
and products of the images of the elements of $M$ are twisted by a $2$-cocycle $M \times M \rightarrow \kx$.
 Both \eqref{AAR} and \eqref{Rsmash} are specializations of \eqref{PxR}. An alternative description of $R$ is as an iterated Ore extension over $\CC(R)$ of the form
$$
R = \CC(R) [x_{l_1}; \sig'_{l_1}] [x_{l_2}; \sig'_{l_2}] \cdots [x_{l_s}; \sig'_{l_s}],
$$
where $F_x(R) = \{ l_1 < l_2 < \cdots < l_s \}$.

\subsection{Main theorem on unipotent automorphisms}
\label{4.2}
Recall that a \emph{connected graded algebra} is a nonnegatively graded algebra $R= \bigoplus_{n=0}^\infty R^n$ 
such that $R^0 = \KK$. For such an algebra,
set $R^{\geq m} := \bigoplus_{n = m}^\infty R^n$ for all $m \in \Znn$. We have used the notation $R^n$ for homogeneous components to avoid conflict with the notation $R_k$ for partial iterated Ore extensions (\S \ref{genCGL}). The algebra $R$ is called \emph{locally finite} if all of its homogeneous components $R^d$ are finite dimensional over $\KK$.

Suppose $R$ is a CGL extension as in \deref{CGL}. Every group homomorphism 
$$
\pi : \xh \to \Zset
$$
gives rise to an algebra $\Zset$-grading on $R$, such that $u \in R^{\pi(\chi_u)}$ for all $\HH$-eigenvectors $u$ in $R$. This makes the algebra
$R$ connected graded if and only if $\pi(\chi_{x_j}) > 0$ for all $j \in [1,N]$. A homomorphism with this property exists 
if and only if the convex hull of $\chi_{x_1}, \ldots, \chi_{x_N}$ in $\xh$ 
does not contain $0$. 
  
\bde{unip} We call an automorphism $\psi$ of 
a connected graded algebra $R$
\emph{unipotent} if 
$$ 
\psi(x) - x \in R^{\geq m+1}, \quad
\forall x \in R^m, \; m \in \Znn.
$$
It is obvious that those automorphisms form a subgroup of $\Aut(R)$, which will be denoted 
by $\UAut(R)$.
\ede

\bth{uaut}
Let $R$ be a symmetric saturated CGL extension which is a connected graded algebra 
via a homomorphism $\pi : \xh \to \Zset$. Then the restriction of every unipotent 
automorphism of $R$ to the core $\CC(R)$ is the identity. 

In other words, every unipotent automorphism $\psi$ of $R$ satisfies
\begin{align}
\label{pr1}
&\psi(x_k) = x_k, \quad &&\forall k \in C_x(R), 
\\ 
&\psi(x_i) = x_i + a_i, \quad &&\forall i \in F_x(R),
\label{pr2}
\end{align}
where for every $i \in F_x(R)$, $a_i$ is a normal element of $R$ lying in $R^{\geq \deg x_i +1}$  
such that $a_i x_i^{-1}$ is a central element of $R[E(R)^{-1}]$.
\eth

The proof of \thref{uaut} is given in \S \ref{4.3}.

The restriction of a unipotent automorphism to $\AA(R)$ can have a very general form
as illustrated by the next two remarks. 

\bre{1}Consider the quantum affine space algebra 
$$
R = \Oq(\KK^3) := \KK \lcor x_1,x_2, x_3  \mid x_i x_j = q x_j x_i \; \forall i < j \rcor,
$$
for a non-root of unity $q \in \kx$,
which is a symmetric CGL extension with respect to the 
natural action of $(\kx)^3$. In this case, $\AA(R) = R$ and $\CC(R) = \KK$. All the generators $x_i$ are prime, 
thus $P_x(R) = \{ 1,2,3 \}$. Introduce the 
grading such that $x_1$, $x_2$, $x_3$ all have degree $1$. The unipotent automorphisms of this algebra are determined \cite[Th\'eor\`eme 1.4.6]{AC} by
$$
\psi(x_1) = x_1, \; \;  
\psi(x_2) = x_2 + \xi x_1 x_3, \; \;
\psi(x_3) = x_3, \quad \text{for some} \; \xi \in \KK.
$$
In particular, in this case the normal element $a_2 = \xi x_1 x_3$ is generally nonzero.
At the same time, the normal elements $a_1$ and $a_3$ vanish.
\ere

\bre{2} It is easy to see that the polynomial algebra $R = \KK[x_1, \ldots, x_N]$ is 
a symmetric CGL extension with the standard action of $(\kx)^N$. In this case,
again we have $\AA(R)= R$. Currently, little is known 
for the very large group of unipotent automorphisms of the polynomial algebras 
in at least 3 variables. 
\ere

In \S \ref{graut} we show how one can explicitly describe the full automorphism 
groups of many symmetric saturated CGL extensions $R$ with small factors $\AA(R)$
using \thref{uaut} together with graded methods. These ``essentially noncommutative'' CGL extensions are very rigid; typically, all automorphisms are graded with respect to the grading of \thref{uaut}, and often there are few or no graded automorphisms beyond the diagonal ones. 
These types of CGL extensions are very common in 
the theory of quantum groups.
We illustrate this by giving  a second proof of 
the Launois--Lenagan conjecture \cite{LL1} that states that 
$$
\Aut( \OO_q(M_n(\KK)) \cong \Zset_2 \ltimes (\kx)^{2n -1}
$$
for all $n>1$, base fields $\KK$, and non-roots of unity $q \in \kx$.
Here, the nontrivial element of $\Zset_2$ acts by the transpose 
automorphism $(X_{lm} \mt X_{ml})$ and the torus acts by rescaling the $X_{lm}$. 
This conjecture was proved for $n=2$ in \cite{AC}, for $n=3$ in \cite{LL2} and for all 
$n$ in \cite{Y-LL}. We reexamine this in \S \ref{graut}, reprove it in a new way, 
and give a very general approach to such relationships based on \thref{uaut}.
\medskip

For a simple Lie algebra $\g$, the algebra $\UU_q^+(\g)$ is the subalgebra 
of $\UU_q(\g)$ generated by all positive Chevalley generators 
$E_\al$, $\al \in \Pi$, recall the setting of \exref{qsc}.
The Andruskiewitsch--Dumas conjecture \cite{AD} 
predicted an explicit description of the full automorphism 
group of $\UU_q^+(g)$. This conjecture was proved in \cite{Y-AD}
in full generality. The key part of the conjecture was to show that 
\begin{equation}
\label{UAutUq}
\UAut(\UU_q^+(\g)) = \{ \id \} 
\end{equation}
for the $\Znn$-grading given by $\deg E_\al =1$, $\alpha \in \Pi$. The next proposition 
establishes that $\CC(\UU^+_q(\g)) = \UU_q^+(\g)$ for all 
simple Lie algebras $\g \neq \sl_2$, and thus 
Eq. \eqref{UAutUq} also follows from \thref{uaut}. Since the pieces of the 
proof in \cite{Y-AD} were embedded in the different steps of the 
proof of \thref{uaut}, this does not give an independent second proof 
of the Andruskiewitsch--Dumas conjecture. However, it illustrates 
the broad range of applications of \thref{uaut} which cover the 
previous conjectures on automorphism groups in this area. 

\bpr{Uq+} For all finite dimensional simple Lie algebras $\g \neq \sl_2$, 
base fields $\KK$ and non-roots of unity $q \in\kx$, 
$$
F_x(\UU_q^+(\g)) = \varnothing, \quad
\mbox{i.e.,} 
\quad
\CC(\UU_q^+(\g))= \UU_q^+(\g).
$$
\epr 

\begin{proof} In the setting of \exref{qsc}, the algebra $\UU_q^+(\g)$ coincides with the algebra $\UU^+[w_0]$ for the longest element $w_0$ of the Weyl group of $\g$. Fix a reduced decomposition 
$w_0 = s_{\al_1} \ldots s_{\al_N}$
for $\al_1, \ldots, \al_N \in \Pi$. Define the roots 
$$
\be_1 := \al_1, \; \be_2 := s_{\al_1} (\al_1), \; \ldots, \; \be_N:=s_{\al_1} \ldots s_{\al_{N-1}}(\al_N)
$$
and Lusztig's root vectors 
$$
E_{\be_1} := E_{\al_1}, \; E_{\be_2} := T_{\al_1} (E_{\al_1}), \; 
\ldots, \; E_{\be_N}:=T_{\al_1} \ldots T_{\al_{N-1}}(E_{\al_N})
$$
in terms of Lusztig's braid group action \cite[\S 8.14]{Ja} on $\UU_q(\g)$. The algebra 
$\UU_q^+(\g)$ has a torsionfree CGL extension presentation 
of the form 
$$
\UU_q^+(\g) = \KK [E_{\be_1}] [E_{\be_2}; \sig_2, \de_2] 
\ldots [E_{\be_N}; \sig_N, \de_N]
$$
for some automorphisms $\sig_\bullet$ and skew derivations $\de_\bullet$, 
the exact form of which will not play a role in the present proof.
Since $\be_1, \ldots, \be_N$ is a list of all positive roots of $\g$, 
for each $\al \in \Pi$ there exists $k(\al) \in [1,N]$ such that 
$$
\be_{k(\al)} = \al.
$$ 
By \cite[Proposition 8.20]{Ja},
\begin{equation}
\label{albe}
E_{\be_{k(\al)}} = E_\al, \quad \forall \al \in \Pi.
\end{equation}
Given $\al \in \Pi$, choose $\al' \in \Pi$ which is connected to $\alpha$ in the Dynkin 
graph of $\g$. (This is the only place we use that $\g \neq \sl_2$.)
The Serre relations imply that $E_\al E_{\al'} \neq \xi E_{\al'} E_\al$
for all $\xi \in \KK$. By \eqref{albe},
$$
k(\al) \notin P_x(\UU_q^+(\g)), \quad 
\forall \al \in \Pi. 
$$
Thus, $E_\al \in \CC(\UU_q^+(\g))$ for all $\al \in \Pi$. Since
$\UU_q^+(\g)$ is generated by $\{E_\al \mid \al \in \Pi \}$, we obtain that
$$
\CC(\UU_q^+(\g)) = \UU_q^+(\g).
$$ 
The decomposition \eqref{decomp} then implies that $F_x(\UU^+_q(\g))$ is empty.
\end{proof}

\subsection{Full automorphism groups}
\label{graut}
There is a large class of quantum nilpotent algebras $R$ for which \thref{uaut} applies and $\CC(R) = R$. 
For such $R$, the only unipotent automorphism is the identity. This lack of unipotent automorphisms often combines with other relations to imply that all automorphisms of $R$ are homogeneous with respect to the grading from \thref{uaut}. We flesh out this statement and analyze several examples in this subsection.

\bde{graut}  {\rm
Let $\psi$ be an automorphism of a connected graded algebra $R$. The \emph{degree zero component} of $\psi$ is the linear map $\psi_0 : R \rightarrow R$ such that
$$
\psi_0(x) \; \text{is the degree} \; d \; \text{component of} \; x, \quad \forall x \in R^d, \; \; d \in \Znn .
$$
The automorphism $\psi$ is said to be \emph{graded} (or \emph{homogeneous of degree zero}) if $\psi = \psi_0$, that is, $\psi(R^d) = R^d$ for all $d \in \Znn$.
 }\ede

\ble{trgruni}
Let $R$ be a locally finite connected graded algebra, $\psi$ an automorphism of $R$, and $\psi_0$ the degree zero component of $\psi$. Assume that $\psi(R^d) \subseteq R^{\ge d}$, for all $d \in \Znn$. Then $\psi_0$ is a graded automorphism of $R$, and the automorphism $\psi_0^{-1} \psi$ is unipotent.
\ele

\begin{proof} It follows immediately from the hypotheses that $\psi_0$ is an algebra endomorphism of $R$. We show that it is an automorphism by proving that $\psi_0$ maps $R^d$ isomorphically onto $R^d$, for all $d \in \Znn$. It suffices to show that $\psi_0(R^d) = R^d$, since $R^d$ is finite dimensional.

Obviously $\psi_0(R^0) = R^0$.  Now assume, for some $d \in \Zset_{>0}$, that $\psi_0(R^j) = R^j$ for all $j \in [0,d-1]$.  Our hypotheses imply that $R^{\ge d} \subseteq \psi^{-1}(R^{\ge d})$, and we next show that this is an equality. If $x \in R \backslash R^{\ge d}$, then $x = y+z$ with $y$ nonzero, $y \in R^j$, and $z \in R^{\ge j+1}$, for some $j \in [0,d-1]$. The assumption $\psi_0(R^j) = R^j$ implies $R^j \cap \ker\psi_0 = 0$, so $\psi_0(y) \ne 0$. Since $\psi(x) - \psi_0(y) \in R^{\ge j+1}$, it follows that $\psi(x) \notin R^{\ge d}$. This shows that, indeed, $R^{\ge d} = \psi^{-1}(R^{\ge d})$, whence $\psi(R^{\ge d}) = R^{\ge d}$. Consequently, any $v \in R^{d}$ can be expressed as $v = \psi(u)$ for some $u \in R^{\ge d}$, and thus $v = \psi_0(u_{d})$ where $u_{d}$ is the degree $d$ component of $u$. This verifies $\psi_0(R^{d}) = R^{d}$ and establishes the required inductive step.

Therefore $\phi_0$ is an automorphism of $R$. It is clear that $\psi_0^{-1} \psi$ is unipotent.
\end{proof}

The condition on $\psi$ in \leref{trgruni} is often satisfied in quantum algebras. In particular, Launois and Lenagan established it in \cite[Proposition 4.2]{LL1} when $R$ is a locally finite connected graded domain, generated in degree $1$ by elements $x_1,\dots, x_n$ such that for all $i \in [1,n]$, there exist $x'_i \in R$ with
$x_i x'_i = q_i x'_i x_i$ for some $q_i \in \KK^*$, $q_i \neq 1$. 
If, in addition, $R$ is a symmetric saturated CGL extension such that $\CC(R) = R$ and $R$ is connected graded via a homomorphism $\pi : \xh \rightarrow \Zset$, we can conclude from \thref{uaut} that all automorphisms of $R$ are graded. We illustrate this by giving a second proof of the descriptions of $\Aut(\Oq(M_{t,n}(\KK)))$ in \cite[Theorem 4.9, Corollary 4.11]{LL1} and \cite[Theorem 3.2]{Y-LL}.

\bpr{qmatr} For all integers $n,t \geq 2$, 
base fields $\KK$, and non-roots of unity $q \in\kx$, 
\begin{equation}
F_x(\OO_q(M_{t,n}(\KK))) = \varnothing, \quad
\mbox{i.e.,} 
\quad
\CC(\OO_q(M_{t,n}(\KK)))= \OO_q(M_{t,n}(\KK)).
\label{qmatrFx}
\end{equation}
Consequently,
$$
\UAut( \OO_q(M_{t,n}(\KK) )) = \{ \id \}
$$
for the grading of $\OO_q(M_{t,n}(\KK))$ with $\deg X_{lm}=1$ for all $l,m \in [1,n]$.
\epr 

\begin{proof} Recall the CGL extension presentation of $R = \OO_q(M_{t,n}(\KK))$
from \eqref{OqMtnCGL} and the function $\eta$ from \exref{OqMtn2}. We have already noted that $R$ is a symmetric CGL extension. The scalars $\la_{kl}$ are all equal 
to powers of $q$. Thus, $R$ is a torsionfree CGL extension, 
and in particular it is saturated.

The only level sets of $\eta$ of cardinality $1$ 
are $\eta^{-1}(n-1)$ and $\eta^{-1}(1-t)$, i.e., the only generators 
of $R$ that are prime are $X_{1n}$ and $X_{t1}$. Thus, $P_x(R) = \{ n, (t-1)n+1 \}$.
The identities 
\begin{align*}
&X_{1,n-1} X_{2n} - X_{2n} X_{1,n-1} = (q-q^{-1}) X_{1n} X_{2,n-1}
\quad \mbox{and}
\\
&X_{t-1,1} X_{t2} - X_{t2} X_{t-1,1} = (q-q^{-1}) X_{t-1,2} X_{t1}
\end{align*}
imply \eqref{qmatrFx}. The final conclusion of the proposition now follows from \thref{uaut}.
\end{proof}

\bth{Autqmatr}
{\rm [Launois--Lenagan, Yakimov]}
For all integers $n,t \ge 2$, base fields $\KK$, and non-roots of unity $q \in \kx$, 
$$
\Aut(\Oq(M_{t,n}(\KK)) = 
\begin{cases}
\DAut(\Oq(M_{t,n}(\KK)) \cong (\kx)^{t + n - 1},
&\mbox{if $n\neq t$},
\\
\DAut(\Oq(M_{t,n}(\KK)) {\cdot} \{ \id, \tau \} \cong (\kx)^{t + n - 1} \ltimes \Zset_2
&\mbox{if $n=t$},
\end{cases}
$$
where $\tau$ is the transpose automorphism of $\OO_q(M_{n,n}(\KK))$ given by 
$\tau(X_{ij}) = X_{ji}$, for all $i,j \in [1,n]$.
\eth

\noindent{\bf Remark}. In the cases where $t$ or $n$ is $1$, $\Oq(M_{t,n}(\KK))$ 
is a quantum affine space algebra. In these cases a description of the automorphism 
groups was found much earlier in \cite{AC} using direct arguments.

\begin{proof}
Let $R = \Oq(M_{t,n}(\KK))$ as in Examples \ref{eOqMtn}, \ref{eOqMtn2}, with $n,t \ge 2$. This algebra is a locally finite connected graded domain in which all generators $X_{ij}$ have degree $1$. By \cite[Corollary 4.3]{LL1}, all automorphisms $\psi$ of $R$ satisfy $\psi(R^d) \subseteq R^{\ge d}$, for all $d \in \Znn$. Thus, by \leref{trgruni}, \thref{uaut}, and \prref{qmatr}, all automorphisms of $R$ are graded. 

It remains to show that any graded automorphism $\psi$ of $R$ has the stated form. We first look at the induced automorphism $\ol{\psi}$ on the abelianization $\ol{R} := R/[R,R]$. Note that the cosets in $\ol{R}$ of the generators $X_{ij}$ satisfy 
$$
\ol{X}_{ij} \ol{X}_{lm} = 0 \; \; \text{if} \; \begin{cases} (i=l, \; j \ne m), \; \text{or}\\  (i \ne l, \; j=m), \; \text{or}\\  (i < l, \; j > m), \end{cases}
$$
and that the $\ol{X}_{ij}^{\,2}$ together with the products $\ol{X}_{ij} \ol{X}_{lm}$ for $i < l$ and $j < m$ form a basis for $\ol{R}^{\, 2}$. It is easily checked that the degree $1$ part of the annihilator of $\ol{X}_{1n}$ has dimension $tn-1$, as does that of $\ol{X}_{t1}$, while no degree $1$ elements of $\ol{R}$ other than scalar multiples of $\ol{X}_{1n}$ or $\ol{X}_{t1}$ have this property. Thus, $\ol{\psi}(\ol{X}_{1n})$ must be a scalar multiple of either $\ol{X}_{1n}$ or $\ol{X}_{t1}$, and similarly for $\ol{\psi}(\ol{X}_{t1})$. It follows that in $R$, we have $\psi(X_{1n}), \psi(X_{t1}) \in \kx X_{1n} \cup \kx X_{t1}$.

Now define
$$
C_s(x) := \{ y \in R^1 \mid xy = q^s yx \}, \; \; \forall s \in \Zset, \; x \in R^1,
$$
and observe that $\psi(C_s(x)) = C_s(\psi(x))$. Since $C_1(X_{1n})$ and $C_1(X_{t1})$ have dimensions $t-1$ and $n-1$, respectively, we conclude that 
\begin{equation}
\label{diag}
\psi(X_{1n}) \in \kx X_{1n} \quad \mbox{and} \quad \psi(X_{t1}) \in \kx X_{t1}
\end{equation}
if $t \neq n$. If $t = n$ and $\psi(X_{1n}) \in \kx X_{t1}$, 
$\psi(X_{t1}) \in \kx X_{1n}$, then the composition $\psi \tau$ will have the property 
\eqref{diag}. Thus, it remains to show that every graded automorphism $\psi$
of $R$ that satisfies \eqref{diag} is a diagonal automorphism. 
It follows from \eqref{diag} that $\psi$ preserves the space 
$$
V := R^1 \cap C_{-1}(X_{1n}) = \KK X_{11} +\cdots+ \KK X_{1,n-1}.
$$
For $j \in [1,n-1]$, the elements $v \in V$ for which $\dim_\KK(V \cap C_1(v)) = n - j - 1$ and $\dim_\KK (V \cap C_{-1}(v)) = j-1$ are just the nonzero scalar multiples of $X_{1j}$. Hence, $\psi(X_{1j}) \in \kx X_{1j}$ for all $j \in [1,n]$. Similarly, $\psi(X_{i1}) \in \kx X_{i1}$ for all $i \in [1,t]$.

Finally, for $i \in [2,t]$ and $j \in [2,n]$, the elements of $C_1(X_{i1} )\cap C_1(X_{1j})$ are exactly the nonzero scalar multiples of $X_{ij}$. We conclude that $\psi(X_{ij}) \in \kx X_{ij}$ for all $i \in [1,t]$, $j \in [1,n]$, showing that $\psi$ is a diagonal automorphism of $R$.
\end{proof}

We give two additional examples which can be established in similar fashion, leaving details to the reader.

\bex{qspo}
First, let $R := \Oq(\mathfrak{sp}\, k^{2n})$ be the quantized coordinate ring of $2n$-dim\-en\-sional symplectic space, with generators $x_1,\dots,x_{2n}$ and relations as in \cite[\S 1.1]{M}. (This presentation gives a symmetric CGL extension presentation, whereas the original presentation in  \cite[Definition 14]{RTF} is not symmetric.) Then:

For all integers $n>0$, base fields $\KK$, and non-roots of unity $q \in \kx$,
$$
\Aut( \Oq(\mathfrak{sp}\, k^{2n})) = \DAut( \Oq(\mathfrak{sp}\, k^{2n}) \cong (\kx)^{n+1}.
$$

Now let $R := \Oq(\mathfrak{o}\, k^m)$ be the quantized coordinate ring of $m$-dim\-en\-sional euclidean space, with generators $x_1,\dots,x_m$ and relations as in \cite[\S\S 2.1, 2.2]{M}. Then:

For all integers $n>0$, base fields $\KK$, and non-roots of unity $q \in \kx$,
\begin{align*}
\Aut( \Oq(\mathfrak{o}\, k^{2n})) &= \DAut( \Oq(\mathfrak{o}\, k^{2n})) {\cdot} \lcor \tau \rcor \cong (\kx)^{n+1} \rtimes \Zset_2,  \\
\Aut( \Oq(\mathfrak{o}\, k^{2n+1})) &= \DAut( \Oq(\mathfrak{o}\, k^{2n+1})) \cong (\kx)^{n+1},
\end{align*}
where $\tau$ is the automorphism of $\Oq(\mathfrak{o}\, k^{2n})$ that interchanges $x_n$, $x_{n+1}$ and fixes $x_i$ for all $i \ne n,n+1$.
\eex

\subsection{Proof of \thref{uaut}}
\label{4.3}
The proof of \thref{uaut} is based on the rigidity of quantum tori 
result of \cite{Y-AD}. This proof is carried out in 
six steps via Lemmas \ref{l1}--\ref{l6} below.
Some parts of it are similar to the proof of the Andruskiewitsch--Dumas 
conjecture in \cite[Theorem 1.3]{Y-AD}, other parts are different.
Throughout the proof we use the general facts for 
CGL extensions established in \cite{GY,GY2}. 

Note that the $\Znn$-grading on the algebra $R$ in \thref{uaut} extends to a $\Zset$-grading on $R[E(R)^{-1}]$, since $E(R)$ is generated by homogeneous elements.
 
\ble{1} In the setting of Theorem {\rm\ref{tuaut}}, for every $k \in [1,N]$
there exists $z_k \in Z(R[E(R)^{-1}])^{\geq 1}$ such that 
$$
\psi(x_k) = (1 + z_k) x_k.
$$
\ele

\noindent {\bf Note}. It follows from \prref{center} that the elements $z_k$ satisfy
\begin{equation}
\label{zk}
z_k \in \NN(R)[E(R)^{-1}], \; \; \forall k \in [1,N].
\end{equation}

\begin{proof} Fix $k \in [1,N]$. There exists an element $\tau$ 
of the subset $\Xi_N$ of the symmetric group $S_N$ defined in 
\eqref{tau} such that $\tau(1)= k$. For example,  one can choose 
$$
\tau =[k, k+1, \ldots, n, k-1, k-2, \ldots, 1]
$$
in the one-line notation for permutations. For the corresponding sequence of 
prime elements, we have $y_{\tau,1} = x_k$. 
The corresponding embeddings $\AA_\tau \subseteq R \subset \TT_\tau$ 
are $X(\HH)$-graded. We use the homomorphism $\pi : \xh \to \Zset$ 
to obtain a $\Znn$-grading on $\AA_\tau$ and a $\Zset$-grading on $\TT_\tau$
for which all generators $y_{\tau, 1}, \ldots, y_{\tau, N}$ have positive degree. 
The embeddings $\AA_\tau \subseteq R \subset \TT_\tau$ become 
$\Zset$-graded. It follows from \prref{sat} that $\TT_\tau$ is a saturated quantum torus 
since $R$ is a saturated CGL extension.

Applying the 
rigidity of quantum tori result in \cite[Theorem 1.2]{Y-AD} and the conversion 
result \cite[Proposition 3.3]{Y-LL}, we obtain that
$$
\psi(y_{\tau, k}) = (1+ c_k) y_{\tau, k} \; \;
\mbox{for some} \; \; c_k \in Z(\TT_\tau)^{\ge 1}, \; \; \forall k \in [1,N].
$$   
By \prref{center}, 
$Z(\TT_\tau) = Z(R[E(R)^{-1}])$. Using that $y_{\tau,1}= x_k$ and setting
$z_k:= c_1$ leads to the desired result.
\end{proof} 

From now on, all characters will be computed with respect to the 
the torus $\DAut(R)$, recall \deref{daut} and the discussion after it.
For an algebra $A$, we denote by $A^*$ its group of units.
 
\ble{2} In the setting of Theorem {\rm\ref{tuaut}}, the elements $z_k \in Z(R[E(R)^{-1}])$,
$k \in[1,N]$, from Lemma {\rm\ref{l1}} define a group homomorphism 
$$
X(\DAut(R)) \to \Fract(Z(R[E(R)^{-1}]))^*
$$
such that
\begin{equation}
\label{gr-hom}
\chi_{x_k} \mt 1 +z_k
\end{equation}
for $k \in [1,N]$.

Consequently, if $u$ is any homogeneous 
element of $R$ and $\chi_u = j_1 \chi_{x_1} + \cdots + j_N \chi_{x_N}$
for some $j_1, \ldots, j_N \in \Zset$, then    
\begin{equation}
\label{psiy}
\psi(u) = (1+z_1)^{j_1} \cdots (1+z_N)^{j_N} u. 
\end{equation}
\ele

\begin{proof} It follows from \cite[Theorem 5.5]{GY} that $X(\DAut(R))$ is generated by 
$\chi_{x_1}, \ldots, \chi_{x_N}$. For $l \in [1,N]$, denote by $X(\DAut(R))_l$ 
the subgroup of $X(\DAut(R))$ generated by $\chi_{x_1}, \ldots, \chi_{x_l}$.

We show by induction on $l$ that there exists a group 
homomorphism $X(\DAut(R))_l \to \Fract(Z(R[E(R)^{-1}]))^*$ satisfying Eq. \eqref{gr-hom} for $k \in [1,l]$. The statement is obvious 
for $l=1$. Assume its validity for $l-1$, where $l \ge 2$. If $\de_l = 0$, then by 
\cite[Theorem 5.5]{GY},
$$
X(\DAut(R))_l = X(\DAut(R))_{l-1} \oplus \Zset \chi_{x_l}
$$
and the statement follows trivially. Now consider the case $\de_l \neq 0$. 
Choose $j < l$ such that $\de_l(x_j) \neq 0$, in other words, $Q_{lj} \neq 0$. 
Choose a monomial $x_{j+1}^{m_{j+1}} \ldots x_{l-1}^{m_{l-1}}$ 
which appears with a nonzero coefficient in $Q_{lj}$, and observe that
$$
\chi_{x_l} = - \chi_{x_j} + m_{j+1} \chi_{x_{j+1}} + \cdots + m_{l-1} \chi_{x_{l-1}}.
$$
The inductive step thus amounts to proving that
\begin{equation}
(1+ z_l) = (1+ z_j)^{-1}(1+ z_{j+1})^{m_{j+1}} \ldots (1+ z_{l-1})^{m_{l-1}}.
\label{toprove}
\end{equation}
The inductive hypothesis, the fact that $z_1, \ldots, z_l$ are central, and the fact that all monomials appearing with nonzero coefficients in $Q_{lj}$ have the same $X(\DAut(R))$-degrees give
$$
\psi(Q_{lj})=
(1+ z_{j+1})^{m_{j+1}} \ldots (1+ z_{l-1})^{m_{l-1}} Q_{lj}.
$$
Applying $\psi$ to the identity $Q_{lj} = x_l x_j - \la_{lj} x_j x_l$ and again using 
that $z_1, \ldots, z_l$ are central leads to
\begin{align*}
(1+ z_{j+1})^{m_{j+1}} \ldots (1+ z_{l-1})^{m_{l-1}} Q_{lj} &=
(1+z_l)(1+z_j) (x_l x_j - \la_{lj} x_j x_l) 
\\
&= (1+z_l)(1+z_j) Q_{lj}. 
\end{align*}
This implies \eqref{toprove} because $Q_{lj} \neq 0$, and completes the induction, establishing the first part of the lemma.

The last statement of the lemma follows from the first part of the lemma, the centrality of the $z_k$, and the fact that all monomials $x_1^{m_1} \cdots x_N^{m_N}$ appearing with nonzero coefficients in $u$ have the same $X(\DAut(R))$-degree as $u$.
\end{proof} 

\ble{3} Any symmetric CGL extensions $R$ of length $N$ is a free left $\NN(R)$-module in which $\NN(R) x_k$ is a direct summand, for all $k \in [1,N] \backslash P_x(R)$. 

If $u x_k \in R$ for some $u \in \NN(R) [E(R)^{-1}]$ and $k \in [1,N] \backslash P_x(R)$, 
then $u \in \NN(R)$.
\ele 

\begin{proof} Theorem 4.11 in \cite{GY} proves that $R$ is a free left module over $\NN(R)$ and constructs an 
explicit basis of it. For $k \in [1,N] \backslash P_x(R)$, the element $x_k$ becomes 
one of the basis elements, because $|\eta^{-1}(\eta(k))| > 1$. This proves the first part of the lemma. 

For the second part, write $r := ux_k$ and $u = e^{-1}y$ for some $e \in E(R)$ and $y \in \NN(R)$. Then $er = yx_k \in \NN(R) x_k$. It follows from the first part of the lemma that $r \in \NN(R) x_k$, and therefore $u \in \NN(R)$.  
\end{proof}
 
\ble{4} In the setting of Theorem {\rm\ref{tuaut}}, the elements $z_k$
from Lemma {\rm\ref{l1}} satisfy
$$
z_k \in Z(R)^{\geq 1}, \quad \forall k \in [1,N] \backslash P_x(R).
$$
\ele 

\begin{proof} By \eqref{zk}, $z_k \in \NN(R)[E(R)^{-1}]$. Furthermore, 
$z_k x_k = \psi(x_k) - x_k \in R$. We apply the second part of \leref{3} to $u := z_k$ to obtain $z_k \in \NN(R)$ and so 
$$
z_k \in R \cap Z(R[E(R)^{-1}])^{\geq 1} = Z(R)^{\geq 1}
$$
for all $k \in [1,N] \backslash P_x(R)$.
\end{proof}

\ble{5} In the setting of Theorem {\rm\ref{tuaut}}, the elements $z_k$
from Lemma {\rm\ref{l1}} satisfy
$$
z_k = 0, \; \; \forall k \in [1,N] \backslash  P_x(R).
$$
\ele

\begin{proof} Let $k \in [1,N] \backslash P_x(R)$ and denote
$$
\eta^{-1} (\eta(k)) = \{ k_1 < \ldots < k_m \}.
$$
By \thref{CGL}, $y_{k_m}$ is a homogeneous 
prime element of $R$ and 
$$
\chi_{y_m} = \chi_{x_{k_1}} + \cdots + \chi_{x_{k_m}}.
$$
Applying \eqref{psiy} with $u= y_{k_m}$ gives
$$
\psi(y_{k_m}) = (1+z_{k_1}) \ldots (1+z_{k_m}) y_{k_m}.
$$
From \leref{4}, $z_{k_1}, \ldots, z_{k_m} \in Z(R)$. 
So, $\psi(R y_{k_m}) \subseteq R y_{k_m}$. At the same time, $R y_{k_m}$ is a height
one prime ideal of $R$, and and so $\psi(R y_{k_m})$ is a 
height one prime ideal. Therefore $\psi(R y_{k_m}) = R y_{k_m}$,
which implies that $(1+z_{k_1}) \ldots (1+ z_{k_m})$ is a unit of $R$. 
The group of units of a CGL extension is reduced to scalars, thus
$$
(1 + z_{k_1}) \ldots (1+ z_{k_m}) \in \kx.
$$
Since $z_{k_1}, \ldots, z_{k_m} \in R^{\geq 1}$, this is only possible 
if $z_{k_1}= \ldots = z_{k_m} = 0$. Therefore $z_k =0$.
\end{proof}
 
\ble{6} In the setting of Theorem {\rm\ref{tuaut}}, the elements $z_k \in Z(R[E(R)^{-1}]$
from Lemma {\rm\ref{l1}} satisfy 
$$
z_k =0, \quad \forall k \in C_x(R).
$$
\ele  

\begin{proof} The statement was proved for $k \in [1,N] \backslash P_x(R)$ in \leref{5}. Now let 
$k \in C_x(R) \cap P_x(R)$. There exist $j, l \in [1,N] \backslash P_x(R)$ 
such that $j < k < l$ and there is a monomial $x_{j+1}^{m+1} \ldots x_{l-1}^{m_{l-1}}$ 
with $m_k >0$ that appears with a nonzero coefficient in $Q_{lj}$. Applying 
$\psi$ to the identity $Q_{lj} = x_l x_j - \la_{lj} x_j x_l$ 
and using Lemmas \ref{l2} and \ref{l5} gives 
$$
(1+z_{j+1})^{m_{j+1}} \ldots (1+ z_{l-1})^{m_{l-1}} = (1 + z_l)(1 + z_j) = 1. 
$$
Since $z_{j+1}, \ldots, z_{l-1} \in Z(R[E(R)^{-1}])^{\geq 1}$ and $R[E(R)^{-1}]$ 
is a graded domain, $z_t=0$ for all $t \in [j+1, l-1]$ 
such that $m_t>0$. Thus $z_k=0$, because $m_k >0$. 
\end{proof}
\medskip

\noindent
\emph{Proof of Theorem {\rm\ref{tuaut}.}} This follows from Lemmas \ref{l1} and \ref{l6}, setting $a_i = z_i x_i$ for $i \in F_x(R)$, recalling that $x_i$ is normal in $R$ for all $i \in F_x(R)$.
\qed


\end{document}